\documentclass[a4paper,latexcad,twoside,11pt]{article}
\usepackage{amssymb} % \mathfrak, \mathbb
\usepackage{graphicx}
\usepackage{epstopdf}
\begin{document}
\title{ New presentations of a link and virtual link }
\author{
Liangxia Wan \thanks{\it E-mail  address: $lxwan@bjtu.edu.cn$. }
 \\
  \small\it Department of Mathematics,
Beijing Jiaotong University, Beijing $100044$, China}

\date{}
\maketitle

\noindent{\small {\bf Abstract} Unique representations of a link and a virtual link are introduced and algebraic systems on links and virtual links are constructed respectively.
Based on the algebraic systems, Reduction Crossing Algorithms for them are proposed which are used to reduce the number of crossings in a link and virtual link by applying a main tool--a pass replacement. For an infinite class unknots $\cal U$, One can transform  each $K$ into a trivial knot in at most  $O(n^c)$ by applying the corresponding algorithm for certain set $\cal U$ of unknots where $c$ is a constant, $K\in {\cal U}$ and $n=|V(K)|$. As special consequences, three unknots are unknotted  which are Goeritz's unknot,Thistlethwaite's unknot and Haken's unknot (image courtesy of Cameron Gordon). Moreover, an infinite family of unknots $K_{G_{2k,2l}}\in {\cal U}$ are unknotted in $O(nloglogn)$ time.

%\noindent MSC(2000): 05C10, 05C30}

\vskip 5mm
\noindent {\bf $1.$ Introduction}
\vskip 5mm

\noindent

 A link consists of $m$ closed non-self-intersecting curves embedded in $R
^3$ for a positive integer $m$. If $m=1$, then it is called a knot. A link projection is obtained by projecting a link to $R^2$. If there are a finite double crossings (no any other crossings) in a projection and the over-crossing line and the under-crossing line transverse at each crossing in the plane, then the projection is called a {\it diagram }of a link. In fact, this presentation of a link can be dated to Brunn's work \cite{Br98}. There are several notations such as Dowker notation \cite{DT83}, Gauss code and Conway notation \cite{Co70} to represent a knot or a link with certain combinatorial or algebraic structures. However, two Conway notations may represent the same underlying knot \cite{MMR08} and one Dowker notation (or one Gauss code \cite{Ga00})may represent two distinct knots \cite{Ad94}(See Fig.1).

This paper introduces a new presentation based on a diagram of a link and a planar embedding of a graph. Recall that the {\it rotation} $\sigma_u$ at the vertex $u$ is a cyclic permutation of its incident edges. Then ${\cal P}_G=\prod\limits_{u\in V(G)} \sigma_u$ is a {\it rotation system }of
$G$. Edmonds found that there is a bijection between
the rotations systems of a graph and its combinatorial orientable embeddings \cite{Ed60}. Youngs
provided the first proof published \cite{Yo63}. The idea of the
bijection above can be dated to some of the earlier work about
imbeddings such as Dyck \cite{Dy88} and Heffter \cite{He98}.

Given a diagram  $L$ of a link, if one labels every crossing with two mutual letters $a$ and $a^-$ where $a$ and $a^-$ represent the overcrossing and undercrossing respectively, then one obtains a marked diagram which is still denoted by $L$.  Its shadow  is a planar embedding of a marked 4-regular graph obtained from $L$ by regarding $a$ and $a^-$ as the same vertex $a$(no overcrossings or undercrossings at crossings). Here, its incident edges at $a$ are $e_i=(a^r,x_i^{r_i})$ such that $e_1e_3,e_3e_1\notin \sigma_a$ where $i=1,3$ for $r=+$ and there is not any other crossing between $a^r$ and $x_i^{r_i}$ along the corresponding curves of $L$ for $r,r_i\in \{+,-\}$ with $1\le i\le 4$.
Since topologically planar embeddings of a graph are determined by ${\cal P}_{G}$ and the {\it infinite face} $f_{G}$ which is the unbounded face, the marked planar embedding with an infinite face is called its {\it embedding} representation. Through this paper a link $L$, unless otherwise indicated, is always such an embedding representation $L$ that each vertex has an anticlockwise rotation, although it would be more precise to call it a marked diagram representative of a link in $R^3$.
%\vskip 3mm

Obviously, given a marked diagram $L$ of a link, there exists one and only one embedding presentation. Conversely, without loss of generalization, suppose that $a$ is a crossing of an embedding representation and suppose further that its incident edges are $e_i=  (a,x_i^{r_i})$ with $i=1,3$ and $e_i=(a^-,x_i^{r_i})$ with $i=2,4$ for $r_i\in\{+,-\}$ where $e_1e_3,e_3e_1\notin \sigma_a$. Let $e_1$ and $e_3$ on the same line and let the crossing be an overcrossing $a$, and
simultaneously let $e_2$ and $e_4$ on the same line and let the crossing be an undercrossing $a^-$. Then one and only one diagram is obtained. Thus we arrive at the following result.
\vskip 3mm
\noindent{\bf Theorem $1.1.$} {\it There are a bijection between marked diagrams and embedding representations.}

\vskip 3mm
The isomorphism of links is reasonable according to the isomorphism of topological embedding.
\vskip 3mm
\noindent{\bf Definition $1.2.$} {\it For two links $L_1$ and $L_2$, $L_1\cong L_2$ if and only if there exists a bijection $\phi:V(L_1)\rightarrow V(L_2)$ such that $\phi(a^-)=(\phi(a))^-$, $\phi(\sigma_{a(L_1)})=\sigma_{\phi(a)(L_2)}$ for every $a\in V(L_1)$,  and $\phi(f_{L_1})=f_{L_2}$ where $f_{L_i}$ is the infinite face of $L_i$ for $1\le i\le 2$.}

\vskip 2mm
Let ${\cal L}$ be the set of links.
Now we construct a new algebraic system $({\cal L},\sim)$ where an equivalence relation $"\sim"$ is defined blow.
\vskip 3mm
{\bf $\Omega0$:}
Suppose that $e_x=(x_1^{r_1},y_1^{s_1})$  and $e_{y}=(x_2^{r_2},y_2^{s_2})$ are on the same face of $L$ with $r_j,s_j\in \{+,-\}$ for $1\le j\le 2$. A new link  for $s\in \{+,-\}$
$$L^{+x^sy^s}\sim L$$
where $L^{+x^sy^s}$ follows from $L$ by adding vertices $x,y$ and then replacing edges $e_x$ and $e_y$ with their subdivisions.
Here, if $e_x\ne e_y$, then $(x_1^{r_1},x^s)$, $(x^s,y^s)$, $(y^s,y_1^{s_1})$ is the subdivision of $e_x$ and $(y_2^{s_2},x^{-s})$, $(x^{-s},y^{-s})$, $(y^{-s},x_2^{r_2})$ is the subdivision of $e_y$. Otherwise,  $(x_1^{r_1},x^s)$, $(x^s,y^s)$, $(y^s,y^{-s})$, $(y^{-s},x^{-s})$, $(x^{-s},y_1^{s_1})$ is the subdivision of $e_x$.
\vskip 3mm
$\Omega1$:
Set $U=\{x,y\}$ and
set $(x,x^-),(y,y^-)\notin E(L)$ for $L\in {\cal L}$.  If $e_+,e_-\in E(L)$ where $e_r=(x^r,y^r)$ for $r\in \{+,-\}$, then
$$ L^{-U}\sim L.$$
Suppose that other incident edges of $x$ and $y$ are $(x^r,z_r^{s_r})$ and  $(y^r,v_r^{t_r})$ for   $r,s_r,t_r\in \{+,-\}$ respectively.
$L^{-U}$ is obtained from $L$ by deleting vertices $x$, $y$ and their incident edges, adding an edge $(z_r^{s_r},v_r^{t_r})$ with $z_r^{s_r}\neq y^r$ for $r\in \{+,-\}$, adding two edges $(x,x^-)_i$ with $z_+^{s_+}= y$ and then adding two edges $(y,y^-)_i$ with $z_-^{s_-}= y^{-}$ for $1\le i\le 2$.

\vskip 2mm
{\it Case $1.$} $f_L\ne (e_+,e_-)$;
\vskip 2mm
{\it Case $2.$} otherwise.

\vskip 3mm
{\bf $\Omega2$:}
 Set $e=(x,x^-)\in E(L)$ for $L\in {\cal L}$. Suppose that other incident edges of $x$ are $(x^r,y_r^{s_r})$  for $r,s_r\in \{+,-\}$. If $y_+^{s_+}\ne x^-$, then
$$L^{-x}\sim L$$
where $L^{-x}$ follows from $L$ by deleting the crossing $x$ and its incident edges and then adding the edge $(y_+^{s_+},y_-^{s_-})$.
\vskip 2mm
{\it Case $1$.} $f_L\ne (e)$;

{\it Case $2$.} otherwise.
\vskip 3mm
{\bf $\Omega3$:}
If a triangle $\triangle=((x^{r},y^{r})$, $(y^{-r},z^{-s})$, $(z^{s},x^{-r}))$ is a face of a link $L$ where $r,s\in \{+,-\}$, then $$L^\triangle\sim L$$
where  ${\cal P}_{L^\triangle}={\cal P}_{L}^{\{U;W\}}$ and $f_{ L^\triangle}=f_{ L}^{\{U;W\}}$. Here, $U=\{x^{r},y^{r},y^{-r},z^{-s},z^{s},x^{-r}\}$ and
$W=\{y^{r},x^{r},z^{-s},$

\noindent $y^{-r},x^{-r},z^{s}\}$.

\vskip 2mm
{\it Case $1.$} $f_L\ne\Delta$;
\vskip 1mm
{\it Case $2.$} otherwise.
\vskip 2mm

We state the following results.
\vskip 3mm
 \noindent{\bf Theorem $1.3.$} {\it  If $L$ is connected, then $L$ is uniquely determined by ${\cal P}_{L}$.}
\vskip 3mm
Since each diagram $L$ of a link in $R^3$ is the union of some connected diagrams $L_i$ for $l\ge 2$ and $1\le i\le l$ such that $f_L=\bigcup\limits_{1\le i\le l}f_{L_i}$ by using a generalized exchange move \cite{Dy06}, we actually verify the following conclusion.
\vskip 3mm
\noindent{\bf Theorem $1.4.$} {\it  A link $L$ is uniquely determined by ${\cal P}_{L}$.}
\vskip 3mm
The following result is concluded for the equivalence relation $"\sim"$ on the algebraic system $({\cal L},\sim)$.
\vskip 3mm
\noindent{\bf Theorem $1.5.$} {\it  For any $L_1,L_2\in {\cal L}$,  $L_1\sim L_2$ if and only if $L_1$ can transform into $L_2$ by a sequence of $\Omega i$ for $0\le i\le 3$.}
\vskip 2mm
 We also consider the unknot recognition problem in this paper. There has been much progress since  Haken provided the first algorithm by introducing a normal surface \cite{Ha61} in 1961. In 1999 Hass, Lagarias and Pippenger showed that Unknot recognition is in NP \cite{HLP99}.  In 2003 Dynnikov gave an algorithm to unknot an arc rectangular of an unknot by monotonic (sometime not \cite{La15}) simplifications \cite{Dy06}. Lackenby provided an algorithm with a polynomial upper bound on Reidemeister moves by combining Dynnikov's methods with the use
of normal surfaces in 2015 \cite{La15}. Despite the progress, it is still tough to unlink unknots such as Goeritz's unknot \cite{Go34},Thistlethwaite's unknot and Haken's unknot (image courtesy of Cameron Gordon) illustrated by Lackenby in \cite{La15}. The main reason is that its binding weight of a normal compression disc used in Lackenby's algorithm is generally exponential.

In this paper we directly perform operations on links, find a pass replacement as a powerful tool and then provide  the Reduction Crossing Algorithm. Set ${\cal U}$ to be a set of such unknots that each $K\in {\cal U}$ is equivalent to an $K'\in {\cal U}$ for $|\mu(K)|\ge 3$. Here, $K'$ is obtained by applications of equal pass replacements  in at most $O(n^{c_1})$ time  and there exists one of a short pass replacement and $\Omega i$  in  $K'$  where $c_1$ is a constant independent of $n$.
 Then $K$ is transformed into a trivial knot in $O(n^c)$ time for any $K\in {\cal U}$ where $c=\max\{n^3loglogn,c_1+1\}$. %$c=\max\limits_{K\in {\cal U}}\{n^3loglogn,c_1+1\}$.
 We arrive at $c=2$ for a special infinite class of unknots $K_{G_{2k,2l}}$ with $k>0,l\ge 0$. Moreover, we transform three known unknots above into a trivial knot. Similarly, we introduce the embedding presentations of a virtual link, an oriented link and an oriented virtual link.

This paper is organized as follows. In Section $2$, some related preliminaries are reviewed and some notations are explained. In Section $3$, we study the embedding representation of a link and prove theorems $1.3$ and $1.5$. In Section $4$ the Reduction Crossing Algorithm are presented and some examples are provided. In Section $5$ and $6$ we introduce the embedding presentations of a virtual link, an oriented link and oriented virtual link respectively. Finally some open problems are given in Section $7$. In Appendix, the Haken's unknot is unlinked.

\vskip 5mm
\noindent{\bf $2.$ Preliminaries}

\vskip 5mm
In this section, we  state some notations and review some related information about links \cite{Ad94,Re32}, virtual links \cite{Ka99} and embeddings \cite{Li08}.

\vskip 3mm

\noindent {\it 2.1. Notations}

\vskip 3mm
Let $a^+$ and $a^-$ be mutual inverse letters for a letter $a$. In brevity, $a^+=a$ through this paper.
Set $U$ and $W$ to be a set of letters. Then $|U|$ denotes the number of elements in $U$. $U\setminus W$ denotes the set of letters from $U$ by deleting letters in $W$. A {\it linear sequence} $A$ is a sequence of letters and $(A)$ denotes a rotation of letters. If $ab\in (A)$, then $a$ is {\it adjacent} to $b$.
Suppose that a linear sequence $A=abc\cdots z$. Then $A^{re}=z\cdots cba$ is called the reverse of $A$.
 If a rotation of edges $f=((a_{1}^{r_1},a_{2}^{-r_2}), \cdots, (a_{k-1}^{r_{k-1}},a_{k}^{-r_{k}}),(a_{k}^{r_k},a_{1}^{-r_1}))$  forms a face where $a_i\ne a_j$ for $1\le i\ne j\le k$, then $f$ is denoted by $f=(a_1a_2\cdots a_k)$ in brevity. Especially, $f_L$ denotes the infinite face of a link $L$.

Suppose that $X=(A_1)(A_2)\dots (A_l)$ where $A_i$ are linear sequences for $l\ge 1$ and $1\le i\le l$. Let $(A_1)^{-U}$ denote the rotation from $A_1$ by deleting letters $x,x^-$ for any $x\in U$. Generally, let $X^{-U}=(A_1)^{-U}(A_2)^{-U}\cdots (A_l)^{-U}$.
Let $X^{\{a^r;b_1^{s_1},b_2^{s_2},\cdots, b_l^{s_l}\}}$ denote a product of rotations from $X$ by replacing the letter $a^r$ with a sequence of letters $b_1^{s_1},b_2^{s_2}\cdots b_l^{s_l}$ for $r,s_i\in \{+,-\}$, $l\ge 1$ and $1\le i\le l$. Let $X^{\{e;e_1,e_2\cdots e_l\}}$ denote a rotation from $X$ by replacing the edge $e$ with the sequence of edges $e_1,e_2\cdots e_l$. Suppose that $U_i,W_i$ are sets of letters or edges. Let $X^{\{U_1,\cdots, U_l;W_1,\cdots, W_l\}}=X^{\{U_i;W_i|1\le i\le l\}}$ denote the rotation from $X$ by replacing elements of $U_i$ with these of $W_i$ for $1\le i\le l$.

\vskip 3mm

\noindent {\it 2.2. Planar embeddings }

 \vskip 3mm

 Given a graph $G=(V,E)$ with $|V(G)|\ge k+1$ for $k\ge 1$, a sequence of vertices $P=uu_1u_2\cdots u_lv$ for $l\ge 0$ is called a {\it path} from $u$ to $v$ if $(u_i,u_{i+1})\in E(G)$ for $0\le i\le l$ where $u_0=u$, $u_{l+1}=v$ and $u_i\ne u_j$ for $0\le i\ne j\le l+1$. A permutation of letters $C=(a_1a_2\cdots a_k)$ for $k\ge 1$ is a {\it cycle} if $(a_{i},a_{i+1})\in E(\mu G)$ for $1\le i\le k$ where $a_{k+1}=a_1$  and $a_i\ne a_j$ for $1\le i\ne j\le k$. If there is a path from any $u,v\in V(G)$, then $G$ is {\it connected}. Based on Dijkstra's shortest algorithm \cite{Di59},   Thorup's shortest algorithm gives an  $O(mloglogm)$ time bound for the single source shortest path problem on a connected graph with $m$ edges \cite{Th00}.

 If  $G_1=(V(G)\setminus V_1, E(G)\setminus E_1)$ is connected for any set $V_1$ of vertices with $|V_1|=k-1$, then $G$ is $k$-{\it connected} where $E_1$ is the set of edges incident with vertices of $V_1$. Especially, $G$ is connected for $k=1$. The greatest integer $k$ such that $G$ is $k$-connected
is the {\it connectivity} $\kappa(G)$ of $G$. $G$ is disconnected if there are $u,v\in V(G)$ such that there does not exist any path from $u$ to $v$ in this paper. The maximal connected subgraph of $G$ is called a connected {\it component }of $G$. For a graph $G$ with $|V(G)|\ge 2$, if there exists a vertex $u\in V(G)$ such that the
number of connected components of $G_1=(V(G)\setminus \{u\}, E(G)\setminus E_1)$ is greater than that of $G$ where $E_1$ is the set of edges incident with $u$, then $u$ is a {\it cutvertex}.

If $G$ can be drawn in the plane such that no two edges meet in a point
other than a common end, then this drawing is called a {\it planar embedding} $\mu(G)$ of $G$. The {\it infinite face} is the unbounded face and other faces are called {\it interior faces}.
Given two faces $f_i$ of $\mu(G)$ for $1\le i\le 2$, if $f_1$ and $f_2$ have a common edge, then $f_1$ is adjacent to $f_2$. Regard each face as a vertex, define such two vertices to be adjacent if they are adjacent in $\mu(G)$ and then get the {\it dual} graph of $\mu(G)$.

Given a directed graph $\vec{G}=(V,E)$, if $\vec{G}$ can be drawn in the plane such that no two arcs meet in a point
other than a common end where $E$ is the set of arcs, then this drawing is called an  embedding  of $\vec{G}$ on the plane. The infinite face is the unbounded face and other faces are called {\it interior faces}. Let $\sigma_{u(\mu G)}$ denote the rotation at $u$ in $\mu(\vec{G})$.
\vskip 3mm
\noindent{\bf Theorem $2.1$.} ({\it Jordan Curve Theory}\cite{Jo1887}.) Let $C$ be equivalent to a circle on the plane. Then $R^2-C$ has two components.
\vskip 12mm
%\vskip 2mm

\setlength{\unitlength}{0.97mm}
\begin{center}
\begin{picture}(100,85)
\qbezier(-3,89)(15,77)(1,73)
\qbezier(1,85)(-14,77)(1,73)
\qbezier(2.5,86)(5,88)(7,89)
\put(11,81){\line(1,0){16}}
\put(27,81){\line(-4,1){4.1}}
\put(27,78){\line(-1,0){16}}
\put(11,78){\line(4,-1){4.1}}
\put(31,89){\line(0,-1){16}}

\qbezier(63.5,84)(77,77)(63,73)
\qbezier(67,87)(48,77)(63,73)
\qbezier(61.5,85)(60,86)(59,86)

\put(35,81){\line(1,0){16}}
\put(51,81){\line(-4,1){4.1}}
\put(51,78){\line(-1,0){16}}
\put(35,78){\line(4,-1){4.1}}

\begin{footnotesize}

\put(19,82){{$I^-$}}
\put(19,75){{$I$}}
\put(-21,80){{$RI:$}}

\put(42,82){{$I$}}
\put(42,75){{$I^-$}}

\put(-21,60){{$RII:$}}

\put(19,62){{$II^-$}}
\put(19,54){{$II$}}

\put(42,62){{$II$}}
\put(42,54){{$II^-$}}
\put(-21,36){{$RIII:$}}
\put(26,-1){Fig.$0$:Reidemeister moves}
\end{footnotesize}
\qbezier(-5,66)(10,57)(-5,53)
\qbezier(5,65)(3,64.5)(1,63.5)
\qbezier(-0.8,62.5)(-9,58)(-0.5,55)
\qbezier(0,54.5)(2,54)(5,53)
\put(11,60){\line(1,0){16}}
\put(27,60){\line(-4,1){4.1}}
\put(27,57){\line(-1,0){16}}
\put(11,57){\line(4,-1){4.1}}
\put(30,66){\line(0,-1){16}}
\put(33,66){\line(0,-1){16}}
\put(35,60){\line(1,0){16}}
\put(51,60){\line(-4,1){4.1}}
\put(51,57){\line(-1,0){16}}
\put(35,57){\line(4,-1){4.1}}
\qbezier(57,66)(60,65)(62.5,63.5)
\qbezier(64,62.5)(70,58)(65,55)
\qbezier(63,53)(60,52)(57,53)

\qbezier(67,65)(65,64.5)(63,63)
\qbezier(63,63)(58,58)(62,55)
\qbezier(62,55)(64,54)(67,53)
\put(-5,46){\line(2,-5){7}}
\put(4,46){\line(-2,-5){4}}
\put(-1,33){\line(-2,-5){1.8}}
\qbezier(-7,42)(-6,43)(-4.5,43)
\qbezier(-3.5,43.5)(0,46)(2.6,43.5)
\qbezier(3.5,43)(4.5,43)(6.5,41)
\put(11,37){\line(1,0){16}}
\put(27,37){\line(-4,1){4.1}}
\put(27,37){\line(-4,-1){4.1}}
\put(11,37){\line(4,-1){4.1}}
\put(11,37){\line(4,1){4.1}}
\put(33,46){\line(2,-5){7}}
\put(39,46){\line(-2,-5){2.8}}
\put(35.8,38){\line(-2,-5){3.8}}
\qbezier(30,33)(31,32)(32.5,32)
\qbezier(34,31.5)(36,30)(38,31.5)
\qbezier(39.5,32)(41,32.5)(42,33.5)
\put(50,46){\line(2,-5){4}}
\put(55,34){\line(2,-5){2.2}}
\put(59,46){\line(-2,-5){7}}
\qbezier(48,42)(49,43)(50.5,43)
\qbezier(51.5,43.5)(55,46)(57.6,43.5)
\qbezier(58.5,43)(59.5,43)(61.5,41)

\put(66,37){\line(1,0){16}}
\put(82,37){\line(-4,1){4.1}}
\put(82,37){\line(-4,-1){4.1}}

\put(66,37){\line(4,-1){4.1}}
\put(66,37){\line(4,1){4.1}}

\put(88,46){\line(2,-5){2.8}}
\put(91.5,37.5){\line(2,-5){3.6}}
\put(94,46){\line(-2,-5){7}}
\qbezier(85,33)(86,32)(87.5,32)
\qbezier(89,31.5)(91,30)(93,31.5)
\qbezier(94.5,32)(96,32.5)(97,33.5)
\put(-5,26){\line(2,-5){1.5}}
\put(-3,21){\line(2,-5){5}}

\put(4,26){\line(-2,-5){1.5}}
\put(2,20){\line(-2,-5){2}}
\put(-1,13){\line(-2,-5){1.8}}
\qbezier(-7,20)(0,24)(6.5,19)

\put(11,17){\line(1,0){16}}
\put(27,17){\line(-4,1){4.1}}
\put(27,17){\line(-4,-1){4.1}}

\put(11,17){\line(4,-1){4.1}}
\put(11,17){\line(4,1){4.1}}

\put(33,26){\line(2,-5){5.3}}
\put(38.8,11.5){\line(2,-5){1.5}}
\put(39,26){\line(-2,-5){2.8}}
\put(35.8,18){\line(-2,-5){2.3}}
\put(33.3,11.5){\line(-2,-5){1.5}}
\qbezier(30,13)(36,10)(42,13.5)

\put(50,26){\line(2,-5){1.5}}
\put(52,21){\line(2,-5){2}}

\qbezier(48,20)(55,24)(61.5,19)

\put(59,26){\line(-2,-5){1.5}}
\put(57,20){\line(-2,-5){4.7}}

\put(55,14){\line(2,-5){2.2}}

\put(66,17){\line(1,0){16}}
\put(82,17){\line(-4,1){4.1}}
\put(82,17){\line(-4,-1){4.1}}

\put(66,17){\line(4,-1){4.1}}
\put(66,17){\line(4,1){4.1}}

\put(88,26){\line(2,-5){2.8}}
\put(91.5,17.5){\line(2,-5){1.8}}
\put(93.5,11.6){\line(2,-5){1.5}}
\put(94,26){\line(-2,-5){5.3}}
\put(88.5,11.4){\line(-2,-5){1.5}}
\qbezier(85,13)(91,10)(97,13.5)
\end{picture}
\end{center}

%\vskip 2mm

\vskip 3mm
\noindent {\it 2.3 Links}
\vskip 2mm
In 1926, Reidemeister, independently Alexander and Briggs introduced the Reidemeister moves and proved the following results.
\vskip 3mm
\noindent {\bf Theorem $2.2.$} {\it(See \cite{AB26,Re27}.)
Given two diagrams $L_i$ for $1\le i\le 2$, if $L_2$ can change to $L_1$ by a sequence of  Reidemeister moves, then they are diagrams of the same link in $R^3$.
}

\vskip 3mm
Given a diagram of a link, an  overpass is a subarc of a link goes over at least one crossing but never goes under a crossing and the number of overcrossings on it is called its length. An  underpass is a subarc of a link goes under at least one crossing but never goes over a crossing and the number of undercrossings on it is called its length. A bridge is an overpass that could not get longer and a {\it subway} is an underpass that could not get longer.
\vskip 3mm

\noindent  {\it 2.4. Virtual links}
 \vskip 3mm
 In 1999, Kauffman introduced a  {\it virtual link} $L$ which consists of $m$ closed non-self-intersecting curves embedded in $S\times I$ for a positive number $m$ \cite{Ka99}. A virtual link  is a virtual knot for $m=1$.  A virtual link projection is obtained by projecting a virtual link to $R^2$ such that the set of crossings consists of two types, one of which is the set of classical crossings denoted by $V(\Gamma)$ and the other of which is the set of artificial crossings denoted by $V(\Sigma)$. If neither
three points on the project to the same point and nor
vertex projects to the same point as any other point in a virtual link, then the projection is called its {\it virtual diagram}.

\vskip 5mm
\noindent{\bf $3.$ Embedding presentation of a link }
\vskip 5mm
In this section we study the embedding presentation of a link and introduce operations on links, especially a pass replacement.

Since
each link $L$ is the union of some connected links $L_i$ for $l\ge 2$ and $1\le i\le l$ such that $f_{ L}=\bigcup\limits_{1\le i\le l}f_{L_i}$ by using a generalized exchange move \cite{Dy06}.  As a special case a trivial knot is represented by $O$ where ${\cal P}_O=((o,o^-)_1, (o,o^-)_1, (o,o^-)_2, (o,o^-)_2)$, which is reasonable by Reidemeister Move I. Similarly, a trivial link with $m$ components $L=\bigcup\limits_{1\le i\le l}O_i$ such that $f_{L}=\bigcup\limits_{1\le i\le l}f_{O_i}$ where $O_i$ are trivial knots such that $O_i=((o_i,o_i^-)_1, (o_i,o_i^-)_1,(o_i,o_i^-)_2, (o_i,o_i^-)_2)$ for $1\le i\le l$. Therefore, it is enough to study its each connected component if one wants to study a link. Therefore, we regard a link as a connected link. Next we look at an example.

\vskip 3mm
\noindent{\bf Example $3.1.$} Two knots $L_1$ and $L_2$ are shown in Fig.1 (Fig.2.8 of \cite{Ad94}) where
\begin{eqnarray*}
{\cal P}_{L_1} & = & (e_3,e_4,e_2,e_1) (e_1,e_2,e_6,e_5)(e_5,e_6,e_7,e_3) (e_8,e_9,e_{10},e_7)(e_9,e_{8},e_{12},e_{11}) (e_{11},e_{12},e_4,e_{10}),
\end{eqnarray*}
\begin{eqnarray*}
{\cal P}_{L_2} & = & (e_3,e_4,e_2,e_1) (e_1,e_2,e_6,e_5)(e_5,e_6,e_7,e_3) (e_8,e_7,e_{10},e_9)(e_9,e_{11},e_{12},e_{8}) (e_{11},e_{10},e_4,e_{12}),
\end{eqnarray*}
$f_{L_1}=(e_3,e_7,e_{10},e_4)$ and $f_{L_2}=(e_3,e_7,e_8,e_{12},e_4)$.
Here, $e_1=(a_1^-,a_2),e_2=(a_1,a_2^-),e_3=(a_1,a_3^-),e_4=(a_1^-,a_6)$, $e_5=(a_2^-,a_3),e_6=(a_2,a_3^-)$, $e_7=(a_3,a_4^-),$ $e_8=(a_4,a_5^-),$ $e_9=(a_4^-,a_5))$, $e_{10}=(a_4,a_6^-),$  $e_{11}=(a_5^-,a_6),e_{12}=(a_5,a_6^-).$

$L_1\neq L_2$ because ${\cal P}_{L_1}\neq {\cal P}_{L_2}$.  \hskip 10mm $\Box$

\vskip 3mm

%\vskip 3mm

\setlength{\unitlength}{0.97mm}
\begin{center}
\begin{picture}(90,73)
%%%%%%%%%%%%%%%%%%%%%%L_1%%%%%%%%%
%\qbezier(12,38)(38,24)(16,24)
\put(-22,58){\line(1,0){6}}
\put(-14,58){\line(1,0){13}}
%\put(-14,18){\line(1,0){6}}
%\put(-6,18){\line(1,0){6}}
\put(2,58){\line(1,0){6}}
%\put(6,18){\circle{7}}
\qbezier(-7,59)(-2,70)(1,58)
\qbezier(-15,58)(-10,44)(-7,57)
\qbezier(-15,58)(-5,80)(8,58)
%%%%%%%%%%%%
\put(10,58){\line(1,0){6}}
\put(18,58){\line(1,0){13}}
%\put(18,18){\line(1,0){6}}
%\put(26,18){\line(1,0){6}}
\put(34,58){\line(1,0){6}}
\qbezier(1,58)(5,44)(10,58)
\qbezier(17,58)(21,44)(25,57)
\qbezier(25,59)(28.5,70)(33,58)
\qbezier(-22,58)(9,30)(33,58)
\qbezier(17,58)(29,80)(40,58)
\begin{footnotesize}
\put(-17.5,59){{$a_1$}}
\put(-10.5,59){{$a_2$}}
\put(-3.5,59){{$a_3$}}
\put(14,59){{$a_4$}}
\put(21,59){{$a_5$}}
\put(28.5,59){{$a_6$}}
\put(7,38){{$L_1$}}
\end{footnotesize}
\put(48,58){\line(1,0){24}}
\put(72,58){\line(1,0){6}}
%\put(6,18){\circle{7}}
\qbezier(63,58)(68,70)(71,58)
\qbezier(55,58)(60,44)(63,58)
\qbezier(55,58)(65,80)(78,58)
%%%%%%%%%%%%
\put(80,58){\line(1,0){24}}
\put(88,58){\line(1,0){13}}
%\put(18,18){\line(1,0){6}}
%\put(26,18){\line(1,0){6}}
\put(104,58){\line(1,0){6}}
\qbezier(71,58)(75,44)(80,58)
\qbezier(87,58)(92,44)(95,58)
\qbezier(95,58)(98.5,70)(103,58)
\qbezier(48,58)(79,30)(103,58)
\qbezier(87,58)(99,80)(110,58)
\begin{footnotesize}
\put(52.5,59){{$a_1$}}
\put(59.5,59){{$a_2$}}
\put(66.5,59){{$a_3$}}
\put(84,59){{$a_4$}}
\put(92,59){{$a_5$}}
\put(98.5,59){{$a_6$}}
\put(77,38){{$L_1$}}
\end{footnotesize}

\put(-22,18){\line(1,0){6}}
\put(-14,18){\line(1,0){13}}
\put(2,18){\line(1,0){6}}
\qbezier(-7,19)(-2,30)(1,18)
\qbezier(-15,18)(-10,4)(-7,17)
\qbezier(-15,18)(-5,40)(8,18)
%%%%%%%%%%%%
\put(10,18){\line(1,0){6}}
\put(18,18){\line(1,0){13}}
\put(34,18){\line(1,0){6}}
\qbezier(1,18)(5,4)(10,18)
\qbezier(17,18)(21,30)(25,19)
\qbezier(25,17)(28.5,4)(33,18)
\qbezier(17,18)(29,-1)(40,18)
\qbezier(33,18)(37,30)(42,18)
\qbezier(42,18)(16,-20)(-22,18)
\begin{footnotesize}
\put(-17.5,19){{$a_1$}}
\put(-10.5,19){{$a_2$}}
\put(-3.5,19){{$a_3$}}
\put(14,19){{$a_4$}}
\put(21,19){{$a_5$}}
\put(28.5,19){{$a_6$}}
\put(7,-5){{$L_2$}}
\end{footnotesize}

\put(53,18){\line(1,0){24}}
\put(61,18){\line(1,0){13}}
\put(77,18){\line(1,0){6}}
\qbezier(68,18)(73,30)(77,18)
\qbezier(60,18)(65,4)(68,18)
\qbezier(60,18)(70,40)(83,18)
%%%%%%%%%%%%
\put(85,18){\line(1,0){24}}
\put(93,18){\line(1,0){13}}
\put(109,18){\line(1,0){6}}
\qbezier(77,18)(80,4)(85,18)
\qbezier(92,18)(96,30)(100,18)
\qbezier(100,18)(102.5,4)(108,18)
\qbezier(92,18)(104,-1)(115,18)
\qbezier(108,18)(112,30)(117,18)
\qbezier(117,18)(91,-20)(53,18)
\begin{footnotesize}
\put(57.5,19){{$a_1$}}
\put(64.5,19){{$a_2$}}
\put(71.5,19){{$a_3$}}
\put(89,19){{$a_4$}}
\put(96,19){{$a_5$}}
\put(103.5,19){{$a_6$}}
\put(82,-5){{$L_2$}}
\put(35,-9){{Fig.1: $L_1$ and $L_2$}}
\end{footnotesize}
\end{picture}
\end{center}

\vskip 9mm

\hskip 2mm Now we consider the new algebraic system $({\cal L},\sim)$. Let consider relations between the operations $\Omega i$ and Reidemeister moves for $0\le i\le 3$. The following result is easily verified.

\vskip 3mm
\noindent{\bf Lemma $3.2.$} %3.3
{\it There are the correspondences between $\Omega i$ and Reidemeister Moves for $0\le i\le 3$:
$$\Omega 0 \longleftrightarrow II,$$
$$\mbox{ Case 1 of }\Omega 1 \longleftrightarrow II^-,$$
$$\mbox{ Case 1 of }\Omega 2 \longleftrightarrow I^-\mbox{  for a non-trivial component  }$$
and
$$\mbox{ Case 1 of } \Omega 3 \longleftrightarrow\mbox{ RIII. }$$}
\vskip 3mm
$\Omega i$ are mainly distinct from Reidemeister moves in two sides for $0\le i\le 3$. One is  the deletion of the operation $I$ %since it is obtained by using $\Omega 2$ and then using Case 1 of $\Omega 0$.
and the other is the extension of $I^-$, $II^-$ and $RIII$ which is verified in Corollary $3.10$ later.
Given a link $L$, a rotation of letters $X=(a_1^{r_1}a_2^{r_2}\cdots a_l^{r_l})$ is called a {\it curve} if $(a_i^{r_i},a_{i+1}^{r_{i+1}})\in E(L)$ where $a_{l+1}^{r_{l+1}}=a_{1}^{r_{1}}$ for $r_i\in \{+,-\}$, $l\ge 2$ and $1\le i\le l$. Furthermore, if $(x^r,z_1^{\epsilon}),(z_i^{\epsilon},z_{i+1}^{\epsilon})$ and $(z_k^{\epsilon},y^s)\in E(L)$ for $x\ne y$ and $r,\epsilon,s\in \{+,-\}$, then $P=x^rz_1^{\epsilon}z_2^{\epsilon}\cdots z_k^{\epsilon}y^s$ is a {\it pass} of $L$ with the {\it length} $k$ between $x^r$ and $y^s$ for $k\ge 0$ and $1\le i\le k$. If $r=s=-\epsilon$, then $P$ is {\it maximal}. 　Set $e=(x^r,y^s)\in E(L)$. A new link $L_1$ is obtained by adding vertices $z_i$ on $L$ such that $e_0=(x^r,z_1^{\epsilon}),e_i=(z_i^{\epsilon},z_{i+1}^{\epsilon}),e_k=(z_k^{\epsilon},y^s)\in E(L_1)$ for $r,\epsilon,s\in \{+,-\}$, then $e_0,e_1,\cdots,e_k$ is a {\it subdivision} of $e$. In the following part we define a (special) pass replacement of a link.

\vskip 3mm

\noindent{\bf Definition $3.3.$} %3.4
For $L\in {\cal L}$, let $f=(A_1xyA_2zA_3)$ is a face of $L$ for $x\ne y\ne z$ where $A_i$ are linear sequences for $1\le i\le 3$. Suppose that $e=(x^{r},y^{s})$, $e_i=(z^{(-1)^{i}r_1},a_i^{s_i})\in E(L)$ with $r,s,r_1,s_i\in \{+,-\}$ for $1\le i\le 4$  and suppose further that  $e,e_1,e_4\in f$, $a_4\in yA_2z$ and $a_1\in zA_3A_1x$.  A new link $L_{\{e;z\}}$ called a {\it replacement} of $e$ surrounding $z$ is obtained from $L$
where $ L_{\{e;z\}}$
  comes from $L$ by adding vertices $x_i$ and then replacing $e_j$ with their subdivisions for   $1\le i\le 4$ and $0\le j\le 4$. Here, $(x^{r},x_1^{\epsilon}),(x_1^{\epsilon},x_2^{\epsilon}),(x_2^{\epsilon},x_3^{\epsilon}),
  (x_3^{\epsilon},x_4^{\epsilon}),(x_4^{\epsilon},y^{s})$ is the subdivision of $e=e_0$. If $e_{i}\ne e_{j}$ for $1\le i\ne j\le 4$, then $(z^{(-1)^{i}r_1},x_i^{-\epsilon})$, $(x_i^{-\epsilon},a_{i}^{s_i})$ is the subdivision  of $e_i$, otherwise $(z^{(-1)^{i}r_1},x_i^{-\epsilon})$, $(x_i^{-\epsilon},x_{j}^{-\epsilon})$, $(x_{j}^{-\epsilon},z^{(-1)^{j}r_1})$ is the subdivision of $e_i$.
\vskip 3mm
\noindent {\bf Lemma $3.4.$} %3.5
{\it Let $L_{\{e;z\}}$ be a {\it replacement} of $e$ surrounding $z$ as above. Then $$L_{\{e;z\}}\sim L.$$}
\vskip 3mm
\noindent{\bf Proof.} Because $e_i\in f$ for $0\le i\le 1$, by $\Omega 0$
\begin{equation}
L_1 \sim L
\end{equation}
\noindent where $e_{x_1}=e_0, e_{x_3}=e_1$ and $L_1=L^{+x_1^{\epsilon}x_3^{\epsilon}}$ for some $\epsilon\in \{+,-\}$.

Because $e_{x_2}=(x_3^{\epsilon},y^{s})$ and $e_{x_4}=(a_4^{s_4},z^{r_1})$  are on the same face of $L_1$, by applying $\Omega0$
\begin{equation}
L_2 \sim L_1
\end{equation}
\noindent where $L_2=L_1^{+x_2^{\epsilon}x_4^{\epsilon}}$.

Because $\triangle=((x_3^{\epsilon},x_2^{\epsilon}),(x_2^{-\epsilon},z^{r_1}),(z^{-r_1},x_3^{-\epsilon}))$ is an interior face of $L_2$, by $\Omega3$
\begin{equation}
L_2^\Delta \sim L_2.
\end{equation}
 It is obvious that
\begin{equation}
L_2^\Delta =L_{\{e;z\}}
\end{equation}
Combining (1-4), one obtains the conclusion.\hskip 10mm $\Box$
\vskip 3mm

Generally, let $U=\{x_i,x_i^{-}| 1\le i\le k\}$ for a nonnegative integer $k$. For $L\in {\cal L}$, set $(A_1a_1^{r_1}a_2^{r_2}\cdots a_l^{r_l}A_2)$ to be a curve where $A_i$ are linear sequences for $r_j\in\{+,-\}$, $1\le i\le 2$, $l\ge 2$ and $1\le j\le l$. Then $a_l^{r_l}$ is called the $(l-1)$th successor of $a_1^{r_1}$ denoted by $suc^{l-1}(a_1^{r_1})$ and $a_1^{r_1}$ is called the $(l-1)$th precursor of $a_l^{r_l}$ denoted by $pre^{l-1}(a_l^{r_l})$. Specifically, if $l=2$, then $a_1^{r_1}$ is called the {\it precursor} of $a_2^{r_2}$ and $a_2^{r_2}$ is called the {\it successor} of $a_1^{r_1}$, i.e. $a_1^{r_1}=pre(a_2^{r_2})$ and $a_2^{r_2}=suc(a_1^{r_1})$. Let $P_1=x^rx_1^{\epsilon}x_2^{\epsilon}\cdots x_k^{\epsilon}y^s$ is a pass for $r,\epsilon,s\in \{+,-\}$ and  $x,y\notin U$.
Add a pass $P=x^ry_1^{\epsilon}y_2^{\epsilon}\cdots y_m^{\epsilon}y^s$ on $L$ for $m\ge 1$ such that $C=(Px_k^{\epsilon}\cdots x_1^{\epsilon})$ is a cycle. i.e. add vertices $y_j$ for $1\le j\le m$ and replace $e_{y_j}$ with their subdivisions for $0\le j\le m$ such that $x,e_{y_1}\in f_0$, $y,e_{y_m}\in f_{m}$ and $e_{y_{j-1}},e_{y_{j}}\in f_{j-1}$ for $2\le j\le m$ where $f_i$ are faces of $L$ for $0\le i\le m$. Here, $e_0, e_1,\cdots, e_{m}$ is a subdivision of $e_{y_0}=(x^r,y^s)$ where $e_{i}=(y_{i}^{\epsilon},y_{i+1}^{\epsilon})$ with $y_0^{\epsilon}=x^r$ and $y_{m+1}^{\epsilon}=y^s$ for $0\le i\le m$. If $e_{y_p}\ne e_{y_q}$ for any $1\le p< q\le m$, then $e_{a_j}=(y_j^{-\epsilon},a_j^{r_j}),e_{b_j}=(y_j^{-\epsilon},b_j^{s_j})$ is a subdivision of  $e_{y_j}=(a_j^{r_j},b_j^{s_j})$ for $1\le j\le m$. Otherwise, if $e_{y_{i_1}}=e_{y_{i_2}}=\cdots =e_{y_{i_l}}$ with $1\le i_1<i_2\cdots <i_l\le m$ for $l\ge 2$, then $(a_{i_1}^{r_1},y_{i_1}^{-\epsilon}),(y_{i_1}^{-\epsilon}, y_{i_2}^{-\epsilon}), \cdots,  (y_{i_{l-1}}^{-\epsilon}, y_{i_l}^{-\epsilon}), (y_{i_l}^{-\epsilon},b_{i_1}^{s_1})$ is a subdivision of $e_{y_{i_1}}$. Thus a planar embedding $\mu(G)$ is constructed.
\vskip 3mm

\noindent{\bf Definition $3.5.$}
Let $P_1$ be a pass above. A new link
$L_P$
is called a {\it pass replacement } of $L$ where $P$ is given above. Here, $L_P$ forms from $\mu(G)$ by deleting vertices $x_i$ and their incident edges, adding edges $(pre^p(x_i^{-\epsilon}),suc^q(x_i^{-\epsilon}))$ with $pre^p(x_i^{-\epsilon}),suc^q(x_i^{-\epsilon})\notin U$, $pre^l(x_i^{-\epsilon}), suc^t(x_i^{-\epsilon})\in U$ for each $1\le l\le p-1$, $1\le t\le q-1$ and $1\le i\le k$. Especially, if some curve consists of $x_{i_1}^{-\epsilon}, x_{i_2}^{-\epsilon},\cdots, x_{i_l}^{-\epsilon}$ for $1\le j\le 2$, $l\ge 2$ and $1\le i_1<i_2<\cdots<i_l\le k$, then add a trivial component $O_{i_1}$ in $L_P$.
If $m\le k-1$, then $L_P$ is a {\it short pass replacement} of $L$. If $m=k$, then $L_P$ is an {\it equal pass replacement} of $L$. Otherwise, $L_p$ is a {\it long pass replacement}.

\vskip 3mm
\noindent {\bf Lemma $3.6.$}
 {\it For any $L\in {\cal L}$, let $L_P$ be a pass replacement for any nonnegative integer $k$ and  positive integer $m$ as Definition $3.5$. Then
$$L_P\sim L.$$}
\vskip 3mm
\noindent{\bf Proof.} Consider the planar embedding $\mu(G)$ above.
We verify this conclusion by induction on the number $n$ of vertices contained in the interior of $C$.

Firstly, consider the case $n=0$. Now verify this case by induction on $k$. If $k=0$, then $m$ is even by the Jordon Curve Theorem. Next induction on $m$. For $m=2$, $((y_1^{\epsilon},y_2^{\epsilon}), (y_1^{-\epsilon},y_2^{-\epsilon}))$ is an interior face of $\mu(L_P)$ by using the Jordon Curve Theory. By case 1 of $\Omega1$
$$L_P^{-U}\sim L_P$$
where $U=\{y_j|1\le j\le 2\}.$
Since $L_P^{-U}=L$, the case is clear.

Suppose that the conclusion holds for any even number less than $m$. Next consider the case $m$.
Because the cycle $C$ does not contain any vertex, there exists some $1\le l\le m-1$ such that $(y_l^{\epsilon},y_{l+1}^{\epsilon})$ and $(y_l^{-\epsilon},y_{l+1}^{-\epsilon})$ form an interior face of $L_P$. By applying the case 1 of $\Omega1$ on the link $L_P$
\begin{equation}
L_P^{-U}\sim L_P
\end{equation}
where $U=\{y_l,y_{l+1}\}.$
By the induction assumption
\begin{equation}
L_P^{-U}\sim L.
\end{equation}
Combining (5) and (6), one gets
$$L_P\sim L.$$
Therefore the conclusion is obtained by induction for the case that the cycle $C$ does not contain any vertex for $k=0$.

Assume that the conclusion holds for $n=0$, any positive integer less than $k$ and any positive integer $m$ such that $L_P$ is a pass replacement of $L$. Consider the case $k$. If there exists some $1\le l\le m-1$ (or $1\le l\le k-1$) such that $((y_l^{\epsilon},y_{l+1}^{\epsilon}),(y_l^{-\epsilon},y_{l+1}^{-\epsilon}))$ (or  $((x_l^{\epsilon},x_{l+1}^{\epsilon}),(x_l^{-\epsilon},x_{l+1}^{-\epsilon}))$) is an interior face of $L_P$.
Similarly, by applying (5-6) and the induction assumption, one arrives at the conclusion. Otherwise, because the interior of cycle $C$ doesn't contain any vertex, then $k=m$ and  $(x_i^{-\epsilon},y_i^{-\epsilon})\in E(\mu G)$. Let $U=\{x_i,x_i^{-}|1\le i\le k\}$ and $W=\{y_i,y_i^{-}|1\le i\le k\}$. Therefore
$$
L_P= L
$$
where ${\cal P}_{L_P}={\cal P}_{L}^{\{U,W\}}$ and $f_{L_P}=f_{L}^{\{U,W\}}$.

Secondly, assume that the conclusion holds for any positive number less than $n$.

Thirdly, we verify the case $n$.
There exists a vertex $z$ in the interior of $C$ such that both $z$ and at least one of $(x^r,x_1^{\epsilon})$, $(x_i^{\epsilon},x_{i+1}^{\epsilon})$ and $(x_k^{\epsilon},y^s)$ are in the same face of $\mu(G)$ for $1\le i\le k-1$. We further suppose that $z$ and $e_{x_l}=(x_l^{\epsilon},x_{l+1}^{\epsilon})$ in the same face $f_0$ of $L$ for some $1\le l\le k-1$ and leave other cases to readers to verify in a similar argument. Suppose that  $e_{c_i}=(z^{(-1)^{i}\epsilon_0},c_{i}^{\epsilon_j})\in E(L)$ with some $\epsilon_j\in \{+,-\}$ for $1\le i\le 4$ and $0\le j\le 4$ and suppose further that $f_0=(B_1x_lx_{l+1}B_2zB_3)$, $(z^{-\epsilon_0},c_{1}^{\epsilon_{1}}), (z^{\epsilon_0},c_{4}^{\epsilon_{4}})\in f_0$, $c_4\in x_{l+1}B_2z$ and $c_1\in zB_3B_1x_l$. Without loss of generalization, set $e_{c_i}\ne e_{c_j}$ for $1\le i\ne j\le 4$.
A replacement $L_{\{e_{x_l};z\}}$ of $e_{x_l}$ surrounding $z$ is obtained such that $(x_l^{\epsilon},x_{k+1}^{\epsilon}), (x_{k+1}^{\epsilon}, x_{k+2}^{\epsilon}), (x_{k+2}^{\epsilon},x_{k+3}^{\epsilon}), (x_{k+3}^{\epsilon},x_{k+4}^{\epsilon}),$\\ $(x_{k+4}^{\epsilon},x_{l+1}^{\epsilon})$ is a subdivision of $e_{x_{l}}$ and that
$(z^{(-1)^i{\epsilon_0}},$ $x_{k+i}^{-\epsilon}), (x_{k+i}^{-\epsilon},c_i^{\epsilon_i})$ is a subdivision of $e_{c_i}$ for each $1\le i\le 4$.
Here,
$L_{\{e_{x_l};z\}})$ is constructed from
  $L$ by adding vertices $x_{k+i}$ and then replacing $e_{x_l}$ and $e_{c_i}$ with their subdivisions for $1\le i\le 4$.
By Lemma $3.4$,
\begin{equation}
L_{\{e_{x_l};z\}}\sim L.
\end{equation}

Set $W_3=\{j|e_{c_i}=e_{y_j}\mbox{ for some }1\le i\le 4,1\le j \le m\}$ where each $e_{y_j}$ is defined above for $1\le j\le m$. Now construct a planar embedding $\mu(G_1)=(V,E)$ from $L_{\{e_{x_l};z\}}$ by adding vertices $y_j$, edges $(y_l^{\epsilon},y_{l+1}^{\epsilon})$ for $0\le l\le m$ and other incident edges of $y_j$ for $1\le j\le m$  % $W_3=\{j|(z^{(-1)^{i+1}\epsilon_0},c_i^{\epsilon_i})=(a_j^{r_j},b_j^{s_j})\mbox{ for some }1\le %i\le 4,1\le j \le m\}$, $W_4=\{j|c_i^{\epsilon_i}=b_j^{r_j},z^{(-1)^{i+1}\epsilon_0}=a_j^{s_j}\mbox{ for %some }1\le i\le 4,1\le j \le m\}$.
 such that if $j\in W_3$, then let $(x_{k+i}^{-\epsilon},y_j^{-\epsilon}),  (y_j^{-\epsilon},c_i^{\epsilon_i})$  be a subdivision of
$e_{c_i}$; otherwise let $e_{a_j}$, $e_{b_j}$ is a subdivision of $e_{y_j}$. Here, $y_0^{\epsilon}=x^r$, $y_{m+1}^{\epsilon}=y^s$.

Consider the cycle $C_1=(Px_k^{\epsilon}\cdots x_{l+1}^{\epsilon}x_{k+4}^{\epsilon}\cdots x_{k+1}^{\epsilon}x_{l}^{\epsilon}\cdots x_1^{\epsilon})$ in $\mu(G_1)$. $C_1$ contains $n-1$ crossings. Obviously, $L_P$ is also a pass replacement of $L_{\{e_{x_l};z\}}$ between $x^r$ and $y^s$. By induction assumption
\begin{equation}
L_{P}\sim L_{\{e_{x_l};z\}}.
\end{equation}
Therefore the case follows from (7-8).

Finally, the result is verified by induction on $n$. \hskip 10mm $\Box$

\vskip 2mm

 The following conclusion is easily derived from the case $k=0$ of Lemma $3.6$.
 \vskip 3mm
\noindent{\bf Corollary $3.7.$} %3.8
{\it Given a link $L$, set a pass $P_1=x^{r}y_1^{\epsilon}y_2^{\epsilon}\cdots y_m^{\epsilon}y^{s}$ for $1\le j\le 2$, $r,\epsilon,s\in \{+,-\}$ and a positive integer $m$. If $x$ and $y$ is on the same face, then
a new link
$$L_P\sim L$$
where $P=x^ry^s$.}
\vskip 3mm
The following results consider the relations of two links when one differs from the other only in there infinite faces.
\vskip 3mm
\noindent{\bf Lemma $3.8.$} %3.9
{\it Let $f_1$ be an interior face of a connected link $L$. If $f_1$ is adjacent to $f_L$, then a new link $$L_1\sim L$$
 such that ${\cal P}_{L_1}={\cal P}_{L}$  and $f_{L_1}=f_1.$}
 \vskip 3mm

 \noindent{\bf Proof.} If $f_L$ has only one edge, suppose that $f_L=(e)$ where $e=(x^r,x^{-r})$ for some $r\in \{+,-\}$. Set $e_1=(x^{-r},a_1^{r_1})$ and $e_2=(x^{r},a_2^{r_2})$ to be other incident edges of $x$. If $e_1=e_2$, then $L$ is trivial. Therefore, we consider the case $e_1\ne e_2$.  Let $P=x^rx_1x_2x^{-r}$ such that $e_{3}=(x^r,x_1)$, $e_4=(x_1,x_2)$, $e_{5}=(x_2,x^{-r})$ is the subdivision of $e$, and $e_{6}=(x^{-r},x_1^-)$, $e_{7}=(x_1^-,a_1^{r_1})$ is the subdivision of $e_1$ and such that $e_{8}=(x^{r},x_2^-)$, $e_{9}=(x_2^-,a_2^{r_2})$ is the subdivision of $e_2$.  By Lemma $3.6$
 \begin{equation} L_P\sim L \end{equation}
 where $V(L_P)=V(L)\cup \{x_i|1\le i\le 2\}$, $\sigma_{x(L_P)}=\sigma_{x(L)}^{\{e_1;e_{6}\}}$, $\sigma_{a_i(L_P)}=\sigma_{a_i(L)}^{\{e_1;e_7\}}$, $\sigma_{a_2(L_P)}=\sigma_{a_2(L)}^{\{e_2;e_9\}}$,  $\sigma_{u(L_P)}=\sigma_{u(L)}$ for any $u\in V(L)\setminus \{x,a_1,a_2\}$ and $f_{L_P}=f_1^{\{U_i;W_i|1\le i\le 3\}}$. Here, $U_1=\{e_1\}$, $U_2=\{e\}$, $U_3=\{e_2\}$, $W_1=\{e_{7}\}$, $W_2=\{e_{3},e_{5}\}$ and $W_3=\{e_{9}\}$.

 Because $\triangle=((x_1,x_2),(x_2^-,x^r),(x^{-r},x_1^-))$ is an interior face of $L_P$, by the case 1 of $\Omega3$
\begin{equation} L^\triangle\sim L_P \end{equation}
where $V(L^\triangle)=V(L)\cup \{x_i|1\le i\le 2\}$,  $\sigma_{x(L^\triangle)}=(e_6,e_8,e_1,e_2)$,  $\sigma_{u(L_P)}=\sigma_{u(L)}$ for any $u\in V(\mu L)\setminus \{x\}$ and $f_{L^\triangle}=f_1^{\{U_2;W_4\}}$. Here, $W_4=\{e_6,e_{10},e_4,e_{11},e_8\}$ where $e_{10}=(x_1,x_1^-)$ and $e_{11}=(x_2,x_2^-)$. Though the location of $x_i$ change for $1\le i\le 2$, we still denote an edge in $L^\triangle$ with the same notation as that in $L_P$ if they have the same ends.

Because $(e_{x_i})$ is an interior faces of $L^\triangle$ where $e_{x_i}=(x_i,x_i^-)$ for each $1\le i\le 2$, by the case 1 of $\Omega2$
\begin{equation}L_1\sim L^\triangle\end{equation}
where $L_1=((L^\triangle)^{-x_1})^{-x_2}$, $V(L_1)=V( L)$, ${\cal P}_{L_1}={\cal P}_{L}$
and $f_{L_1}=f_1.$

Thus this conclusion holds for $f_L$ with only one edge by $(9-11)$.

Next suppose that $f_L=(A_1,e_1,A_2,e_2,A_3)$ where $e_i=(x_i^{r_i},y_i^{s_i})$, $e_1\ne e_2$, $A_j$ are sequences of edges for $1\le j\le 3$. Suppose further that $e_1$ is a common edge of $f_1$ and $f_L$.
By $\Omega0$ a new link
\begin{equation}
L_1\sim L
\end{equation}
such that $e_x=(x_1^{r_1},x),e_+=(x,y),e_y=(y,y_1^{s_1})$ is the subdivision of $e_1$ and $e_{y^-}=(x_2^{r_2},y^-),e_-=(y^-,x^-),e_{x^-}=(x^-,y_2^{s_1})$ is the subdivision of $e_2$ where $L_1=L^{+xy}$. Here,  $\sigma_{u(L_1)}=\sigma_{u(L)}$ for any $u\in V(L)\setminus\{x_i,y_i|1\le i\le 2\}$. If $y_2^{s_2}=x_1^{-r_1}$, then $\sigma_{{x_1}(L_1)}=\sigma_{{x_1}(L)}^{\{e_1,e_2;e_x,e_{x^-}\}}$; otherwise $\sigma_{{x_1}(L_1)}=\sigma_{{x_1}(L)}^{\{e_1;e_x\}}$ and $\sigma_{{y_2}(L_1)}=\sigma_{{y_2}(L)}^{\{e_2;e_{x^-}\}}$.
If $y_1^{s_1}=x_2^{-r_2}$, then $\sigma_{{y_1}(L_1)}=\sigma_{{x_2}(L)}^{\{e_1,e_2;e_y,e_{y^-}\}}$; otherwise $\sigma_{{x_2}(L_1)}=\sigma_{{x_2}(L)}^{\{e_2;e_y\}}$ and $\sigma_{{y_1}(L_1)}=\sigma_{{y_1}(L)}^{\{e_1;e_{y^-}\}}$. If $e_2\notin f_1$, then  $f_{ L_1}=f_1^{\{U_5;W_5\}}$, otherwise $f_{L_1}=f_1^{\{U_5,U_6;W_5,W_6\}}$ where $U_5=\{e_1\}$, $U_6=\{e_2\}$, $W_5=\{e_x,e_-,e_y\}$ and $W_6=\{e_{x^-},e_+,e_{y^-}\}.$

Set $U=\{x,y\}$. Because $(e_+,e_-)$ is an interior face of $L_1$, by the case $1 $ of $\Omega1$
\begin{equation}
L_2\sim L_1
\end{equation}
where $L_2=L_1^{-U}$, ${\cal P}_{L_2}={\cal P}_{L}$ and $f_{L_2}=f_1$.

By $(12-13)$, this case holds. Thus this conclusion is right. \hskip 10mm $\Box$

 \vskip 3mm

% \vskip 2mm
\noindent{\bf Lemma $3.9.$} {\it Let $L$ be a connected link. For any interior face $g$ of $L$, then there is a new link $$L_g\sim L$$
 where ${\cal P}_{L_g }={\cal P}_{L }$  and $f_{L_g}=g.$}
 \vskip 3mm
\noindent{\bf Proof.} If $g$ is adjacent to $f_{L}$, then the conclusion holds by Lemma $3.8$; otherwise there exists a path $P=f_Lf_1f_2\cdots f_lg$ for $l\ge 1$ such that $f_i$ is adjacent to $f_{i+1}$ where $f_0=f_L$ and $f_{l+1}=g$ for $0\le i\le l$. By using lemma $3.8$ $l+1$ times, the result is obvious. \hskip 10mm $\Box$
\vskip 3mm
This induces the following case.
\vskip 3mm
\noindent{\bf Corollary $3.10.$} %3.11
{\it Cases $2$ of $\Omega i$  are right for $1\le i\le 3$.}
\vskip 3mm

\noindent{\bf Proof of Theorem $1.3.$} For a connected link $L$, let $L_g$ be a link in Lemma $3.9$. By Lemma $3.9$ $L_g\sim L.$ This means that $L$ and $L_g$ are diagrams of the same link in $R^3$. Therefore the conclusion holds.
\hskip 6mm $\Box$

\vskip 3mm

\noindent{\bf Proof of Theorem $1.5.$} One can get $I^-$ from the combined use $\Omega 0$ and Case 1 of $\Omega 2$. In fact, this result is correct by Lemma $3.2$ and Theorem $2.2$. This extends the scope of Reidemeister moves which is powerful on $({\cal L}, \sim)$. \hskip 10mm $\Box$

\vskip 3mm
Thus Definition $1.2$ and $\Omega3$, still denote them with the same notation, are reduced as follows.

\vskip 3mm

\noindent{\bf Definition $1.2.$} {\it For two links $L_1$ and $L_2$, $L_1\cong L_2$ if and only if there exists a bijection $\phi:V(L_1)\rightarrow V(L_2)$ such that $\phi(a^-)=(\phi(a))^-$ and $\phi(\sigma_{a(L_1)})=\sigma_{\phi(a)(L_2)}$ for every $a\in V(L_1)$.}
\vskip 3mm
\noindent{\bf $\Omega3$:}
If a triangle $\triangle=((x^{r},y^{r})$, $(y^{-r},z^{-s})$, $(z^{s},x^{-r}))$ is a face of $L$ where $r,s\in \{+,-\}$, then $$L^\triangle\sim L$$
where  ${\cal P}_{L^\triangle}={\cal P}_{L}^{\{U;W\}}$. Here, $U=\{x^{r},y^{r},y^{-r},z^{-s},z^{s},x^{-r}\}$ and
$W=\{y^{r},x^{r},z^{-s},y^{-r},x^{-r},z^{s}\}$.

\vskip 5mm
 \noindent{\bf $4.$ Reduction Crossing Algorithm }
\vskip 5mm

In fact, since each diagram of link in $R^3$ can be transformed into an equivalent diagram $L$ by using Reidemeister move I$^-$, $\Omega4$ and  $\Omega5$ \cite{Re32} such that there does not exist any cutvertex in $L$, each $L$ is regarded to be  $2$-connected in this section. We propose an algorithm to reduce the number of crossings for a link as follows.

Given a link $L$, let $P=x^sx_1^{-s}\cdots x_{k}^{-s}y^s$ be a maximal pass for $k\ge 1$ and $U=\{x_i,x_i^-|1\le i\le k\}$. Obviously, $x_1^{-s}\cdots x_{k}^{-s}$ is a bridge (or a subway). Set $G=(V,E)$ where $V(G)=V(L)\setminus \{x_i|1\le i\le k\}$ and

$$E(G)=(E(L)\cup E_1)\setminus \{(x^s,x_{1}^{-s}),(x_{k}^{-s},y^s),(x_i^{-s},x_{i+1}^{-s})|1\le i\le k-1\}.$$
Here

\hskip 10mm $E_1=\{(pre^p(x_i^{s}),suc^q(x_i^{s}))|pre^p(x_i^{s}),suc^q(x_i^{s})\notin U, pre^l(x_i^{s}),suc^t(x_i^{s})\in U $

\hskip 52mm for each $1\le l\le p-1, 1\le t\le q-1, 1\le i\le k\}$.

Let $f_x$, $f_y$ and $f_i$ denote faces obtained by removing the edge $(x^s,x_1^{-s})$, $(x_k^{-s},y^s)$ and $(x_i^{-s},x_{i+1}^{-s})$ for $1\le i\le k-1$, respectively. The dual graph of $G$ is called an {\it adjacent graph} $G_{x^sy^s}$ of $x^sy^s$.

\newpage
\noindent{\bf Reduction Crossing Algorithm}
\vskip 2mm
Set $L$ to be  a $2$-connected non-alternating link with $n$ crossings, $N=0$.
\vskip 1mm
1. If there exists $e=(x,x^-)\in E(L)$ and $y^s\ne x^-$ where $(x,y^s)$ and $(x^-,z^t)$ are other adjacent edges  of $x$  for $s,t\in \{+,-\}$, then let $n=n-1$, $L=L^{-x}$; otherwise go to step 2.
\vskip 1mm
2. If there exist $(x,y),(x^-,y^-)\in E(L)$, then set $n=n-2$, $L=L^{-\{x,y\}}$ and go to step 1; otherwise go to step 3.

\vskip 1mm
3. Suppose that $Q=x^sx_1^{-s}\cdots x_{k}^{-s}y^s$ is a maximal pass for $s\in \{+,-\}$. Construct an {\it adjacent graph} $G_{x^sy^s}$ of $x^sy^s$ and determine the shortest path $P_{xy}$ from $f_x$ to $f_y$ by Thorup's shortest algorithm. Set the length of $P_{xy}$ to be $m$. If $m\le k-1$, then let $n=n+m-k$ and $L=L_{P_{xy}}$ and go to step 1; otherwise step 5;
\vskip 1mm
4. If $m= k$ and $P_{N}\ne P_i$ for $1\le i\le N-1$, %f_xf_{x_1}\cdots f_{x_{k-1}}f_y
then let $L=L_{P_N}$, let $N=N+1$ and then go to step 1; otherwise step 5;
\vskip 1mm
5. (1) If there exists two incident edges between   $x$ and $x^-$ for each $x\in V(L)$,  then stop and $L$ is an unlink;

(2) If $L$ is an alternate link, then stop and $L$ is alternate;

(3) Otherwise, if $N\le N_0$, then let $N=N+1$ and go to step 1, otherwise stop.

(4) If $L$ is disconnected, then $L$ is a splitting link.

\vskip 3mm
\noindent{\bf Theorem} { $4.1.$} {\it  Set ${\cal U}_n=\{K||V(K)|=n\}$ for $1\le n\le 2$. Let $c_1$ be a constant and  for $n\ge 3$ let
$$\begin{array}{ll}
{\cal U}_n=\{K|\mbox{ There exists }\mbox{ equal path replacements in at most } O(n^{c_1}) \mbox{ time such that }K\sim K_1\\
 \hskip 15mm\mbox{ and then there exists  one of }\Omega i\mbox{ and a short  path  replacement such that } \\
 \hskip 15mm K_1\sim K_2\in {\cal U}_l \mbox{ for }|V(K)|=|V(K_1)|=n\ge 3, 1\le i\le 2\mbox{ and some }l<n\}.
 \end{array}
 $$
 Then $K$ is transformed into a trivial knot in at most $O(n^c)$ time for any $K\in {\cal U}$ by using Reduction Crossing Algorithm where ${\cal U}=\bigcup\limits_{n\ge 1}{\cal U}_n$ and $c=\max\{n^3loglogn,c_1+1\}$.}
 \vskip 3mm
 \noindent{\bf Proof.}
 Now choose any knot $K\in {\cal U}$ with $|V(K)|=n$. According the definition of $\cal U$, there exists a knot $K_1\in {\cal U}$ after using equal path  replacements in at most $O(n^{c_1})$ time such that
 \begin{equation}K\sim K_1\end{equation}
 and there exists a short path replacement in $K_1$ such that
 \begin{equation}K_1\sim K_2\end{equation}
where $K_2\in {\cal U}_l$ for $l<n$.

We know that one spends at most $O(n^{c_1})$ time on doing (14). Consider operations on $K_1$. If there exist $\Omega i$ for $1\le i\le 2$, then one can finish it in $O(n)$ time, otherwise one constructs an adjacent graph of a maximal pass $Q$ in $O(n)$ time and then  in at most $ O(nloglogn) $ time on finding the shortest path. Therefore one spends at most $O(n^2loglogn)$ time on doing $(15)$.
Thus one does (14-15) in $O(n^{c_2})$ time where $c_2=\max\{n^2loglogn,c_1\}$. After repeat (14-15) at most $n$ times,                            one transforms $K$ into a trivial knot. Therefore we transform
$K$ into a trivial knot in at most $O(n^c)$ time.\hskip 10mm$\Box$
\vskip 3mm
\noindent{\bf Remark $4.2$} Unknots $K_{G_{2k,2l}}\in {\cal U}_{11+2k+2l}$ are given in Example $4.5$ for $k,l\ge 0$, which implies that $\cal U$ has infinite elements. In addition, if one replace equal path replacements with path replacements in ${\cal U}_n$ for $n\ge 3$, then the result obviously holds.

\vskip 3mm
Next we unknot  Thistlethwaite's unknot, Goeritz's unknot and an infinite family of unknots as examples.
\vskip 12mm
%\vskip 2mm

\setlength{\unitlength}{0.97mm}
\begin{center}
\begin{picture}(100,55)

\put(45,60){\line(1,0){18}}
\qbezier(63,60)(65,55)(63,46)
\qbezier(45,60)(43,55)(45,50)
\put(45,50){\line(1,0){3}}
\put(50,50){\line(1,0){3}}
\qbezier(46,57)(49,55)(49,45)
\qbezier(62,55)(52,55)(53,42)
\put(55,50){\line(1,0){3}}
\put(58,50){\line(0,-1){18}}
\put(63,46){\line(-1,0){4}}
\put(57,46){\line(-1,0){2}}
\put(55,46){\line(0,-1){12}}
\put(53,42){\line(1,0){1.5}}
\put(55.5,42){\line(1,0){2}}
\put(49,45){\line(0,-1){8}}
\put(49,37){\line(1,0){5}}
\put(55.5,37){\line(1,0){2}}
\put(55,34){\line(-1,0){12}}
\put(58,32){\line(-1,0){15}}
\qbezier(43,57)(40,50)(42,30)
\put(42,30){\line(1,0){18}}
\put(59,42){\line(1,0){6}}
\put(59,37){\line(1,0){3}}
\qbezier(62,27)(62,34)(62,37)
\put(65,30){\line(0,1){12}}
\put(65,30){\line(-1,0){2}}
\put(62,27){\line(1,0){6}}
\put(68.5,27){\line(0,1){25}}
\qbezier(64.5,55)(67,55)(68.5,52)
\qbezier(68.5,27)(69,24)(72,23)
\qbezier(69,27)(72,27)(73,23)
\put(71,21){\line(-1,0){28}}
\qbezier(71,21)(73,22)(73,23)
\put(41,21){\line(-1,0){3}}
\qbezier(40,34)(35,30)(38,21)
\qbezier(41,32)(38,28)(42,21)
\qbezier(73.5,22.5)(74,23)(75,20)
\qbezier(75,20)(60,10)(42,21)
\begin{footnotesize}
\put(58,38){{$a_{11}$}}
\put(51.5,35){{$a_{5}$}}
\put(58.5,43){{$a_{14}$}}
\put(51.5,40){{$a_{4}$}}
\put(58.5,47){{$a_{3}$}}
\put(37.5,35){{$a_{6}$}}
\put(41.5,28){{$a_{13}$}}
\put(57.5,28){{$a_{10}$}}
\put(74,23){{$a_{8}$}}
\put(69.5,28){{$a_{9}$}}
\put(41.5,22.5){{$a_{7}$}}
\put(55,51){{$a_{15}$}}
\put(44,48){{$a_{12}$}}
\put(40,58){{$a_{1}$}}
\put(64,56){{$a_{2}$}}
\put(32,5){{Fig.$2$: Thistlethwaite's unknot $K_T$}}
\end{footnotesize}
\end{picture}
\end{center}

\noindent{\bf Example $4.3$.} Thistlethwaite's unknot $K_T$ in Fig.$2$.
Then
\begin{eqnarray*}
{\cal P}_{K_T} & = & (e_1,e_2,e_4,e_3)(e_1,e_7,e_5,e_6)(e_9,e_5,e_{10},e_8) (e_{11},e_{12},e_{8},e_{13})\\
&  &  (e_{11},e_{16},e_{14},e_{15})(e_{18},e_{14},e_{2},e_{17})(e_{19},e_{20},e_{21},e_{17})
 (e_{20},e_{19},e_{23},e_{22})\\
&  & (e_{22},e_{23},e_{6},e_{24}) (e_{24},e_{27},e_{25},e_{26}) (e_{29},e_{15},e_{28},e_{25})
 (e_{3},e_{4},e_{16},e_{30})\\
&  &  (e_{18},e_{21},e_{26},e_{28})
 (e_{9},e_{12},e_{29},e_{27})
  (e_{7},e_{30},e_{13},e_{10})
\end{eqnarray*}
where $e_1=(a_1,a_2),e_2=(a_1^-,a_{6}),e_3=(a_1^-,a_{12}),e_4=(a_1,a_{12}^-),$ $e_5=(a_2,a_{3}^-),e_6=(a_2^-,a_{9}),e_7=(a_2^-,a_{15}),$ $e_8=(a_3^-,a_{4}),e_9=(a_3,a_{14}),e_{10}=(a_3,a_{15}^-),$
$e_{11}=(a_4,a_5),e_{12}=(a_4^-,a_{14}^-),e_{13}=(a_4^-,a_{15}),$
$e_{14}=(a_5,a_{6}^-),e_{15}=(a_5^-,a_{11}^-),e_{16}=(a_5^-,a_{12}),$
$e_{17}=(a_6^-,a_{7}^-),e_{18}=(a_6,a_{13}),$
$e_{19}=(a_7,a_8^-),e_{20}=(a_7^-,a_{8}),e_{21}=(a_7,a_{13}^-),$
$e_{22}=(a_9,a_8^-),e_{23}=(a_9^-,a_{8}),e_{24}=(a_9^-,a_{10})$,
$e_{25}=(a_{10},a_{11}^-),e_{26}=(a_{10}^-,a_{13}),$ $e_{27}=(a_{10}^-,a_{14}^-),$
$e_{28}=(a_{11},a_{13}^-),e_{29}=(a_{11},a_{14}),$ $e_{30}=(a_{15}^-,a_{12}^-).$

Let $P_1=a_{13}a_{16}^-a_{17}^-a_{18}^-a_{15}$ such that $e_{31}, e_{33}$ is a subdivision of $e_{28}$, $e_{35}$, $e_{36}$ is a subdivision of $e_{14}$, $e_{38}$, $e_{39}$ is a subdivision of $e_{16}$.
By Lemma $3.6$ a new knot
$$K_{P_1}\sim K_T$$
where
\begin{eqnarray*}
{\cal P}_{K_{P_1}} & = & (e_1,e_2,e_4,e_3)(e_1,e_7,e_5,e_6)(e_{42},e_{5},e_{10},e_{41})(e_{41},e_{38},e_{35},e_{15})\\
&  &  (e_{18},e_{36},e_{2},e_{17})(e_{19},e_{20},e_{21},e_{17})
 (e_{20},e_{19},e_{23},e_{22})(e_{22},e_{23},e_{6},e_{43})\\
&  &  (e_{42},e_{15},e_{31},e_{43})
 (e_{3},e_{4},e_{39},e_{30})(e_{18},e_{21},e_{32},e_{33})
  (e_{7},e_{30},e_{40},e_{10})\\
&  &  (e_{31},e_{34},e_{33},e_{32})(e_{35},e_{37},e_{36},e_{34})(e_{37},e_{38},e_{40},e_{39}).
\end{eqnarray*}
Here, $e_{31}=(a_{16},$  $a_{11}),e_{32}=(a_{16}^-$, $a_{13}),e_{33}=(a_{16},a_{13}^-)$, $e_{34}=(a_{16}^-,a_{17}^-),$
$e_{35}=(a_5,a_{17})$, $e_{36}=(a_6^-,a_{17})$, $e_{37}=(a_{17}^-,a_{18}^-)$, $e_{38}=(a_5^-,a_{18})$, $e_{39}=(a_{12},a_{18})$, $e_{40}=(a_{15},a_{18}^-)$, $e_{41}=(a_3^-,a_5)$, $e_{42}=(a_3,a_{11})$, $e_{43}=(a_9^-,a_{11}^-)$.

Set  $P_2=a_8a_{19}^-a_{20}^-a_{18}$ such that $(a_1^-,a_{19})$, $(a_{19},a_6)$ is a subdivision of $ (a_{1}^-,a_{6})$ and such that $(a_{17}^-,a_{20})$, $(a_{20},a_{18}^-)$ is a subdivision of $(a_{17}^-,a_{18}^-)$, $P_3=a_1^-a_{21}a_{16}^-$ such that $(a_2^-,a_{21})$, $(a_{21},a_8^-)$ is a subdivision of $(a_{2}^-,a_{8}^-)$, $P_4=a_{15}a_{22}^-a_{23}^-a_{7}$ such that $(a_{3},a_{22})$, $(a_{22},a_{15}^-)$ is a subdivision of $(a_{3},a_{15}^-)$ and such that $(a_{16}^-,a_{23}^-)$, $(a_{23}^-,a_{20})$ is a subdivision of  $(a_{16}^-,a_{20})$, $P_5=a_{23}^-a_{24}a_{15}^-$ such that $(a_{20},a_{24})$, $(a_{24},a_{23})$ is a subdivision of $(a_{20},a_{23})$.

Appling pass replacements in order where $P_i$ are given above  for $2\le i\le 5$ and then appling $\Omega1$ and $\Omega2$, one transforms this unknot into a trivial knot. \hskip 10mm $\Box$
\vskip 3mm
\noindent{\bf Example $4.4$.} Goeritz's unknot $K_{G_{0,0}}$ in Fig.$3$.
Then
\begin{eqnarray*}
{\cal P}_{K_{G_{0,0}}} & = & (e_1,e_4,e_3,e_2)
(e_1,e_7,e_6,e_5)(e_5,e_6,e_9,e_8)
(e_8,e_{11},e_{10},e_{12})\\
&  & (e_{11},e_{b_0},e_{b_0^-},e_{10})
(e_{b_0},e_{9},e_{15},e_{b_0^-})(e_{15},e_{17},e_{16},e_{12})
(e_{17},e_{18},e_{4},e_{16})\\
&  & (e_{18},e_{7},e_{20},e_{19})
(e_{21},e_{22},e_{19},e_{20})(e_{22},e_{21},e_{2},e_{3}).
\end{eqnarray*}
Here $e_1=(a_1,a_2),e_2=(a_1^-,a_{11}),e_3=(a_1,a_{11}^-),e_4=(a_1^-,a_8)$,
$e_5=(a_2^-,a_{3}),e_6=(a_2,a_{3}^-),e_7=(a_2^-,a_9^-)$, $e_8=(a_3^-,a_{4}),e_9=(a_3,a_{6}^-),$
$e_{10}=(a_4,a_{5}^-),e_{11}=(a_4^-,a_5),e_{12}=(a_4^-,a_{7}^-),$ $e_{b_0^-}=(a_5,a_{6}^-),e_{b_0}=(a_5^-,a_6)$, $e_{15}=(a_6,a_7),$ $e_{16}=(a_7,a_{8}^-),e_{17}=(a_7^-,a_8)$, $e_{18}=(a_9,a_8^-),e_{19}=(a_9^-,a_{10}),e_{20}=(a_9,a_{10}^-),$ $e_{21}=(a_{10},a_{11}^-), e_{22}=(a_{10}^-,a_{11})$.

Let a pass $P_1=a_{11}^-y_1a_3^-$ such that $e_{24}$, $e_{25}$ is a subdivision of $e_{18}$. By Lemma $3.6$ a new knot
$$K_{P_1}\sim K_{G_{0,0}}$$ where
\begin{eqnarray*}
{\cal P}_{K_{P_1}} & = & (e_{23},e_{25},e_{26},e_{24})
(e_{27},e_{23},e_9,e_{8})(e_8,
e_{11},e_{10},e_{12})
(e_{11},e_{b_0},e_{b_0^-},e_{10})\\
&  & (e_{b_0},e_{9},e_{15},e_{b_0^-})
(e_{15},e_{17},e_{16},e_{12})(e_{17},e_{24},e_{28},e_{16})
(e_{25},e_{29},
e_{20},e_{19})\\
&  & (e_{21},e_{22},e_{19},e_{20})
(e_{22},e_{21},e_{28},e_{26}).
\end{eqnarray*}
Here, $e_{23}=(y_1,a_3^-),e_{24}=(y_1^-,a_{8}^-),e_{26}=(y_1,a_{11}^-),e_{25}=(y_1^-,a_9)$,
$e_{27}=(a_3,a_9^-)$, $e_{28}=(a_8^-,a_{11})$, $e_{29}=(a_9^-,a_3).$

Let a pass $P_2=a_{8}y_2^-a_5$ such that $e_{30}$, $e_{32}$ is a subdivision of $e_9$. By Lemma $3.6$ an equivalent pass replacement
$$K_{P_2}\sim K_{P_1}$$ where
\begin{eqnarray*}{\cal P}_{K_{P_2}} & = & (e_{30},e_{33},e_{32},e_{31})(e_{23},e_{25},e_{26},e_{24})
(e_{29},e_{23},e_{30},e_{34})
(e_{31},e_{b_0},e_{b_0^-},e_{34})
(e_{b_0},e_{32},e_{35},e_{b_0^-})\\
&  & (e_{33},e_{24},e_{28},e_{35})
(e_{25},e_{27},e_{20},e_{19})
(e_{21},e_{22},e_{19},e_{20})
(e_{22},e_{21},e_{28},e_{26}).
\end{eqnarray*}
Here, $e_{30}=(y_2,a_3),e_{31}=(y_2^-,a_{5}),e_{32}=(y_2,a_{6}^-),e_{33}=(y_2^-,a_8)$,
$e_{34}=(a_3^-,a_{5}^-)$, $e_{35}=(a_6,a_8^-).$

Let a pass $P_3=a_6^-y_3a_9^-$ such that $e_{37}, e_{38}$ is a subdivision of $e_{35}$. By Lemma $3.6$ an equivalent pass replacement
$$K_{P_3}\sim K_{P_2}$$ where
\begin{eqnarray*}
{\cal P}_{K_{P_3}} & = & (e_{39},e_{37},e_{36},e_{38})
(e_{40},e_{25},e_{26},e_{24})(e_{41},e_{b_0},e_{b_0^-},e_{40})
(e_{b_0},e_{36},e_{37},e_{b_0^-})\\
&  & (e_{41},e_{24},e_{28},e_{38})
(e_{25},e_{39},e_{20},e_{19})
(e_{21},e_{22},e_{19},e_{20})
(e_{22},e_{21},e_{28},e_{26}).\\
\end{eqnarray*}
Here, $e_{39}=(y_3,a_9^-)$, $e_{38}=(y_3^-,a_{8}^-)$, $e_{36}=(y_3,a_{6}^-)$, $e_{37}=(y_3^-,a_6)$,
$e_{40}=(y_1,a_5^-)$, $e_{41}=(a_5,a_8).$

Let a pass $P_4=a_6a_9$. By Lemma $3.6$ an equivalent pass replacement
$$K_{P_4}\sim K_{P_3}$$ where
\begin{eqnarray*}
{\cal P}_{K_{P_4}} & = & (e_{42},e_{b_0},e_{b_0^-},e_{43})(e_{b_0},e_{45},e_{44},e_{b_0^-})
(e_{44},e_{45},e_{20},e_{19})
(e_{21},e_{22},e_{19},e_{20})
(e_{22},e_{21},e_{42},e_{43}).
\end{eqnarray*}
Here, $e_{42}=(a_5,a_{11})$, $e_{43}=(a_5^-,a_{11}^-)$, $e_{44}=(a_6,a_9)$, $e_{45}=(a_6^-,a_9^-)$.

By applying $\Omega1$ twice, one gets an trivial knot. \hskip 10mm $\Box$
\vskip 2mm
\vskip 18mm

\setlength{\unitlength}{0.97mm}
\begin{center}
\begin{picture}(100,55)
\qbezier(8,59)(12,65)(16,58)
\qbezier(8,58)(12,52)(16,57)
\put(8,58){\line(-1,1){4}}
\put(17,58){\line(1,1){4}}
\put(7.5,58){\line(-1,-1){7.5}}
\put(16,58){\line(1,-1){8}}
\qbezier(-8.5,50)(-4,43)(-1.5,49.5)
\qbezier(-16,50)(-12,57)(-8.5,50)
\qbezier(-24,50)(-20,57)(-16.5,50)
\qbezier(-8,50)(-4,57)(-0.5,50)
\put(-0.5,50){\line(1,-1){7.5}}
\qbezier(-16.5,50)(-13,43)(-8.5,49)
\qbezier(-24.5,49.5)(-21,43)(-16.5,49)

\qbezier(24.5,50)(28,57)(31.5,50)
\qbezier(32,50)(36,57)(39,51)
\qbezier(25,48.8)(28,43)(32,50)
\qbezier(32,49)(36,43)(39.5,50)

\put(24.5,50){\line(-1,-1){8}}
\qbezier(8,43)(12,49)(16,42)
\qbezier(8,41.5)(12,36)(16.6,42)

\put(7.9,43){\line(-1,-1){5.3}}
\put(16.8,41){\line(1,-1){4}}

\put(4,62){\line(-1,0){28}}
\put(2.6,37.6){\line(-1,0){26.8}}
\put(21,62){\line(1,0){18}}
\put(20.8,36.9){\line(1,0){18.1}}

\qbezier(39,62)(41,57)(39.5,50)
\qbezier(39,36.9)(41,43)(39.5,49)

\qbezier(-24,62)(-28,57)(-24.5,49.5)
\qbezier(-24,37.5)(-28,43)(-24.5,48.5)
\begin{footnotesize}
\put(-28,49.5){{$a_{1}$}}
\put(-21,49.5){{$a_{11}$}}
\put(-7.5,49){{$a_{10}$}}
\put(0.5,49){{$a_{9}$}}
\put(20,49){{$a_{6}$}}
\put(32.5,49){{$a_{5}$}}
\put(39.5,49){{$a_{4}$}}
\put(3.5,57.5){{$a_2$}}
\put(17.5,57){{$a_3$}}
\put(3.5,41.5){{$a_8$}}
\put(17.5,41.5){{$a_7$}}
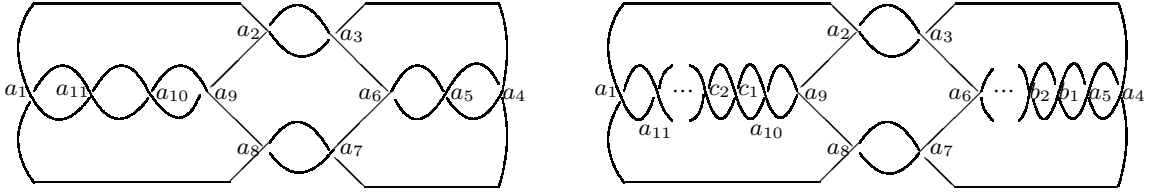
\put(-10,22){{Goeritz's unknot $K_{G_{0,0}}$ }}
\end{footnotesize}

\qbezier(88,59)(92,65)(96,58)
\qbezier(88,58)(92,52)(96,57)
\put(88,58){\line(-1,1){4}}
\put(97,58){\line(1,1){4}}
\put(87.5,58){\line(-1,-1){7.5}}
\put(96,58){\line(1,-1){8}}

\qbezier(75,50)(77,43)(79.5,49.5)
\qbezier(75.5,50)(77,57)(79.5,50)
\put(79.5,50){\line(1,-1){7.5}}

\qbezier(71.5,50.5)(73,57)(75,50)
\qbezier(71,50)(73,43)(75,49)
\qbezier(67.5,50.5)(69,57)(71,50)
\qbezier(67,50)(69,43)(71,49)

\qbezier(56,50)(58,57)(60.5,50)
\qbezier(55.5,49.5)(58,43)(60,49)

\qbezier(60.8,50.5)(61,52.5)(62.5,53.5)
\qbezier(60.5,50)(61,47.5)(62.5,46)

\qbezier(65,46)(66.5,46.5)(66.8,49.5)
\qbezier(65,53.5)(66.5,53)(67,50)

\qbezier(109.5,46)(110.5,46.5)(111,50)
\qbezier(109.5,53.5)(110.5,53)(110.8,50.3)

\qbezier(104.5,50)(105,52)(106,53)
\qbezier(104.5,49)(104.5,47.5)(106,46)

\qbezier(119,50)(121,57)(123,51)
\qbezier(119,49)(121,43)(123.5,50)
\qbezier(115,49)(117,43)(119,50)
\qbezier(115,50)(117,57)(118.5,51)
\qbezier(111,50)(113,57)(114.5,51)
\qbezier(111.5,48.8)(113,43)(115,50)

\put(106.5,50){{\circle*{0.3}}}
\put(107.5,50){{\circle*{0.3}}}
\put(108.5,50){{\circle*{0.3}}}
\put(104.5,50){\line(-1,-1){8.8}}
\qbezier(88,43)(92,49)(96,42)
\qbezier(88,41.5)(92,36)(96,41.5)

\put(87.9,43){\line(-1,-1){5.4}}
\put(96.8,41){\line(1,-1){4}}

\put(84,62){\line(-1,0){28}}
\put(82.6,37.6){\line(-1,0){26.8}}
\put(101,62){\line(1,0){22}}%
\put(100.8,36.9){\line(1,0){22.1}}
\put(65,50){{\circle*{0.3}}}
\put(64,50){{\circle*{0.3}}}
\put(63,50){{\circle*{0.3}}}

\qbezier(123,62)(125,57)(123.5,50)%\qbezier(119,62)(121,57)(119.5,50)
\qbezier(123,37)(125,43)(123.5,49)%7.\qbezier(119,38)(121,43)(119.5,49)

\qbezier(56,62)(52,57)(55.5,49.5)
\qbezier(56,37.5)(52,43)(55.5,48.5)
\begin{footnotesize}
\put(52,49.5){{$a_{1}$}}
\put(58,44){{$a_{11}$}}
\put(67.5,49.5){{$c_{2}$}}
\put(71.5,49.5){{$c_{1}$}}
\put(73.1,44){{$a_{10}$}}%\put(72.5,49){{$a_{10}$}}
\put(80.5,49){{$a_{9}$}}
\put(100,49){{$a_{6}$}}
%\put(105,49){{$b_{2m}$}}
\put(111,49){{$b_{2}$}}
\put(115,49){{$b_{1}$}}
\put(119,49){{$a_{5}$}}
\put(123.5,49){{$a_{4}$}}
\put(83.5,57.5){{$a_2$}}
\put(97.5,57){{$a_3$}}
\put(83.5,41.5){{$a_8$}}
\put(97.5,41.5){{$a_7$}}
\put(86,22){{$K_{G_{2k,2l}}$ }}
\put(28,12){{Fig.3: A type of unknots $K_{G_{2k,2l}}$ }}
\end{footnotesize}
\end{picture}
\end{center}

%\vskip 2mm
\noindent{\bf Example $4.5$.} Unknot a type of unknots $K_{G_{2k,2l}}$ for $k,l\ge 0$ in Fig.$3$.
Then
\begin{eqnarray*}
{\cal P}_{K_G} & = & \prod\limits_{1\le i\le 11}\sigma_{a_i}\prod\limits_{1\le i\le 2k}\sigma_{b_i} \prod\limits_{1\le j\le 2l}\sigma_{c_j}
\end{eqnarray*}
where $\sigma_{a_5}=(e_{11},e_{b_{\theta_k}},e_{b_{\theta_k}^-},e_{10})$,
$\sigma_{a_6}=(e_{b_{2k+{\theta_k}}},e_{9},e_{15},e_{b_{2k+{\theta_k}}^-})$,
$\sigma_{a_{10}}=(e_{c_{\theta_l}^-},e_{c_{\theta_l}},e_{19},e_{20})$, $\sigma_{a_{11}}=(e_{c_{2l+{\theta_l}}},e_{c_{2l+{\theta_l}}^-},$ $e_{2},$ $e_{3})$, $\sigma_{b_{i}}=(e_{b_i},e_{b_{i+1}},e_{b_{i+1}^-},e_{b_i^-})$ for $1\le i\le 2k$, $\sigma_{c_{j}}=(e_{c_j^-},e_{c_{j+1}^-},e_{c_{j+1}},e_{c_j})$ for $1\le j\le 2l$, $\sigma_{a_i}$ and $\sigma_{b_0}$ are the same as those in Example $4.4$ for $1\le i\le 4$ and $7\le i\le 9$. Here, $e_{b_1}=(a_{5}^-,b_{1}),$ $e_{b_1^-}=(a_{5},b_{1}^-),$  $e_{b_{2k+1}}=(a_{6},b_{2k}^-),$ $e_{b_{2k+1}^-}=(a_{6}^-,b_{2k}),$ $e_{b_i}=(b_{i-1}^-,b_{i}),$ $e_{b_i^-}=(b_{i-1},b_{i}^-) $ for $2\le i\le 2k$, $e_{c_1}=(a_{10}^-,c_{1}),$  $e_{c_1^-}=(a_{10},c_{1}^-),$ $e_{c_{2l+1}}=(a_{11},c_{2l}^-),$ $e_{c_{2l+1}^-}=(a_{11}^-,c_{2l}),$  $e_{c_j}=(c_{j-1}^-,c_{j}),$ $e_{c_j^-}=(c_{j-1},c_{j}^-) $ for $2\le j\le 2l$ and for $t\in \{k,l\}$
$$\theta_t=\left\{
\begin{array}{ll}
0, \mbox{ if }t=0;\\
1, \mbox{ otherwise. }
\end{array}
\right.
$$

After we apply the same pass replacements $P_i$ as those in Example $4.4$ for $1\le i\le 4$ in sequence, we get
$$K_{1}\sim K_{G_{2k,2l}}$$
where
 \begin{eqnarray*}
{\cal P}_{K_{1}} & = & \prod\limits_{i\in \{5,6,9,10,11\}}\sigma_{a_i}\prod\limits_{1\le i\le 2k}\sigma_{b_i} \prod\limits_{1\le j\le 2l}\sigma_{c_j}.
\end{eqnarray*}
Here,
$\sigma_{a_5}=(e_{42},e_{b_{\theta_k}},e_{b_{\theta_k}^-},e_{43})$,
$\sigma_{a_6}=(e_{b_{2k+{\theta_k}}},e_{45},e_{44},e_{b_{2k+{\theta_k}}^-})$,
$\sigma_{a_9}=(e_{19},e_{44},e_{45},e_{20})$, $\sigma_{a_{11}}=(e_{43},$ $e_{c_{2l+{\theta_l}}},$ $e_{c_{2l+{\theta_l}}^-},e_{42})$,
$\sigma_{10}$, $\sigma_{b_{i}}$ for $1\le i\le 2k$ and $\sigma_{c_{j}}$ for $1\le j\le 2l$ are given above where $e_{42}=(a_5,a_{11})$, $e_{43}=(a_5^-,a_{11}^-)$, $e_{44}=(a_6,a_9)$ and $e_{45}=(a_6^-,a_9^-)$.

By applying $\Omega1$ $2\min\{k,l\}+2$ times and using $\Omega2$ $2\max\{k,l\}-2\min\{k,l\}$ times, we obtain a trivial knot. \hskip 10mm $\Box$

\vskip 3mm
Since each unknot $K_{G_{2k,2l}}$ has four maximal passes for any $k,l\ge 0$, one can find a short pass in $O(nloglogn)$ time and then conduct pass replacement at most $O(n)$ time. Therefore we transform
$K_{G_{2k,2l}}$ into a trivial knot in $O(nloglogn)$ time.

\vskip 5mm
 \noindent{\bf $5.$ Embedding presentation of a virtual link }
\vskip 5mm

Given a virtual diagram $L$, let $V(\Gamma L)$ and $V(\Sigma L)$ denote the set of classical vertices and virtual vertices respectively. Label each crossing with $x$ and $x^-$. If $x\in V(\Gamma L)$, then $x$ and $x^-$ represent overcrossing and undercrossing respectively,  otherwise they represent twice occurrences of $x$. Similarly,
 its marked shadow  is a planar embedding of a marked 4-regular graph obtained from $L$ by regarding $a$ and $a^-$ as the same vertex $a$(no overcrossings or undercrossings at crossings). Here, its incident edges at $a$ are $e_i=(a^r,x_i^{r_i})$  such that $e_1e_3,e_3e_1\notin \sigma_a$ where $i=1,3$ for $r=+$, $i=2,4$ for $r=-$ and there is not any other crossing between $a^r$ and $x_i^{r_i}$ along the corresponding curves of $L$ for $r,r_i\in \{+,-\}$ with $1\le i\le 4$. Then a virtual diagram $L$ where $R$ is the corresponding marked embedding where the unbounded face  is called the {\it infinite face} of $L$, denoted by $f_L$ and  ${\cal P}_{L}=\prod\limits_{x\in V(L)}\sigma_{x}.$ If $V(\Sigma L)=\emptyset$, then a virtual link is classical. Through this paper a virtual link $L$ is always such a marked planar embedding called an embedding representation.
\vskip 3mm

\noindent{\bf Definition $5.1.$} {\it For two virtual links $L_i$ with $1\le i\le 2$, $L_1\cong L_2$ if and only if there exists a bijection $\phi:V(\Gamma L_1)\rightarrow V(\Gamma L_2)$ and $\phi:\{\alpha|\alpha\in V(\Sigma L_1)\}\rightarrow \{\beta^r|\beta\in V(\Sigma L_2), \mbox{ for some }r\in \{+,-\}\}$ such that $\phi(x^-)=(\phi(x))^-$ and $\phi(\sigma_{x(L_1)})=\sigma_{\phi(x)(L_2)}$ for $x\in V(L_1)$   and $\phi(f_{L_1})=f_{L_2}$.}

\vskip 3mm
Given a virtual diagram $L_1$, let $L_2$ is a link from $L_1$ by exchanging $x$ and $x^-$ for a virtual crossing.
Though  $L_i$  are distinct embedding representations for $1\le i\le 2$, $L_1\cong L_2$ according to the definition $5.1$. From this view point, the following result is given.

\vskip 3mm

\noindent{\bf Theorem $5.2.$} {\it There exists a bijection between marked virtual diagrams and their embedding representations. }
\vskip 3mm
 Some notations without explanation are the same as above.
Set $e=(x^r,y^s)\in E(L)$ for some $r,s\in \{+,-\}$. Delete $e$ and then add edges $(x_i^{\epsilon_i},x_{i+1}^{\epsilon_{i+1}})$ on $L$ with $\epsilon_i\in \{+,-\}$ for $l\ge 1$ and $0\le i\le l$ where $x_0^{\epsilon_0}=x^r$ and $x_{l+1}^{\epsilon_{l+1}}=y^s$. If either $x_j\in V(\Sigma L)$ or $x_j\in V(\Gamma L)$ and $\epsilon_j=\epsilon$ for $1\le j\le l$, then the sequence $(x^{r},x_{1}^{\epsilon_1})$, $(x_1^{\epsilon_1},x_{2}^{\epsilon_2})$, $\cdots,$ $(x_l^{\epsilon_l},y^{s})$ is called a {\it subdivision} of $e$.
Let  $\bar{\cal L}$ denote the set of virtual links. We define an equivalence relation $"\sim"$ satisfying with the following operations on $\bar{\cal L}$:

\vskip 3mm
{\bf $\mho0$:} Suppose that $e_x=(x_1^{r_1},y_1^{s_1})$ and $e_y=(x_2^{r_2},y_2^{s_2})$ are on the same face of $L$ with $r_j,s_j\in \{+,-\}$ for $1\le j\le 2$. A new virtual link  with either $x,y\in V(\Gamma L^{+x^ry^s})$ for $r=s$ and  some $s\in \{+,-\}$ or $x,y\in V(\Sigma L^{+x^ry^s})$
$$L^{+x^ry^s}\sim L$$
where $L^{+x^ry^s}$ follows from $L$ by adding vertices $x,y$ and then replacing edges $e_x$ and $e_y$ with their subdivisions.
Here, if $e_x\ne e_y$, then $(x_1^{r_1},x^r)$, $(x^r,y^s)$, $(y^s,y_1^{s_1})$ is the subdivision of $e_x$ and $(y_2^{s_2},x^{-r})$, $(x^{-r},y^{-s})$, $(y^{-s},x_2^{r_2})$ is the subdivision of $e_y$. Otherwise,  $(x_1^{r_1},x^r)$, $(x^r,y^s)$, $(y^s,y^{-s})$, $(y^{-s},x^{-r})$, $(x^{-r},y_1^{s_1})$ is the subdivision of $e_x$.
\vskip 3mm
$\mho1$: Set $U=\{x,y\}$ and set $(x,x^-),(y,y^-)\notin E(L)$ for $L\in {\bar {\cal L}}$. If $e_r\in E(L)$ for $r\in \{+,-\}$ and if either $x,y\in V(\Gamma L)$ and $s=+$  or $x,y\in V(\Sigma L)$ for $s\in \{+,-\}$ where $e_+=(x,y^s),e_{-}=(x^{-},y^{-s})$, then
$$ L^{-U}\sim L.$$
Suppose that other incident edges of $x$ and $y$ are $(x^r,z_r^{s_r})$ and $(y^s,v_s^{t_s})$ for $r,s,s_r,t_s \in \{+,-\}$ respectively.
$L^{-U}$ is obtained from $L$ by deleting vertices $x$, $y$ and their incident edges, adding an edge $(z_+^{s_+},v_s^{t_s})$ with $z_+^{s_+}\neq y^{s}$, adding an edge $(z_-^{s_-},v_{-s}^{t_{-s}})$ with $z_-^{s_-}\neq y^{-s}$, adding two edges $(x,x^-)_i$ with $z_+^{s_+}= y^{s}$ such that $x\in V(\Gamma L^{-U})$ and then adding two edges $(y,y^-)_i$ with $z_{-}^{s_{-}}= y^{-s}$ for $1\le i\le 2$ such that $y\in V(\Gamma L^{-U})$.

\vskip 3mm
{\it Case $1.$} $f_L\ne (e_+,e_{-})$;
\vskip 2mm
{\it Case $2.$} otherwise.
\vskip 3mm

\vskip 3mm
{\bf $\mho2$:}
 Set $e=(x,x^-)\in E(L)$ for $L\in {\cal L}$. Suppose that other adjacent edges of $x$ are $(x^r,y_r^{s_r})$ for $r,s_r\in \{+,-\}$. If $y_+^{s_+}\ne x^-$, then
$$L^{-x}\sim L$$
where $L^{-x}$ follows from $L$ by deleting the crossing $x$ and its incident edges and then adding the edge $(y_+^{s_+},y_-^{s_-})$.
\vskip 2mm
{\it Case $1$.} $f_L\ne (e)$;

{\it Case $2$.} otherwise.

\vskip 2mm

{\bf $\mho3$:} Let a triangle $\triangle=((x^{s_1},y^{s_2})$, $(y^{-s_2},z^{-s_3})$, $(z^{s_3},x^{-s_1}))$ be a face of $L$ with $s_i\in \{+,-\}$ for $1\le i\le 3$.
If either $x,y,z\in V(\Gamma L)$ for $s_1=s_2$ or $x,y,z\in V(\Sigma L)$, then
$$L^\triangle\sim L.$$
Here, ${\cal P}_{L^\triangle}={\cal P}_{L}^{\{U;W\}}$  and $f_{L^\triangle}=f_{L}^{\{U;W\}}$ where $U=\{x^{s_1},y^{s_2},y^{-s_2},$ $z^{-s_3},z^{s_3},$ $x^{-s_1}\}$,
$W=\{y^{s_2},x^{s_1},z^{-s_3},$ $y^{-s_2},x^{-s_1},z^{s_3}\}$.

\vskip 2mm
{\it Case $1.$} $f_L\ne \triangle$;
\vskip 1mm
{\it Case $2.$} otherwise.
\vskip 2mm

{\bf $\mho4$:}
If a triangle $\Lambda=((\alpha^{r},\beta^{s})$, $(\beta^{-s},a^-)$, $(a,\alpha^{-r}))$ is a face of $L$ with $r,s\in \{+,-\}$, $\alpha,\beta\in V(\Sigma L)$ and $a\in V(\Gamma L)$, then
$$L^\Lambda\sim L$$
Here, ${\cal P}_{L^\Lambda}={\cal P}_{L}^{\{U;W\}}$  and $f_{L^\Lambda}=f_{L}^{\{U;W\}}$ where $U=\{\alpha^{r},\beta^{s},\beta^{-s},a^-,a,$ $\alpha^{-r}\}$ and
$W=\{\beta^{s},\alpha^{r},a^-,\beta^{-s},$ $\alpha^{-r},a\}$.

\vskip 2mm
{\it Case $1.$} $f_L\ne \Lambda$;
\vskip 1mm
{\it Case $2.$} otherwise.

\vskip 2mm

Here, $\mho 0$, the cases 1 of $\mho i$ are just generalized Reidemeister moves for $1\le i\le 4$. The cases 2 of those can be verified later.
\vskip 3mm
\noindent{\bf Definition $5.2.$} For $L\in \bar{\cal L}$, let $f=(A_1xyA_2zA_3)$ is a face of $L$ for $x\ne y\ne z$ where $A_j$ are linear sequences for $1\le j\le 3$. Suppose that $e=(x^{r},y^{s})$, $e_i=(z^{(-1)^{i}r_1},a_i^{s_i})\in E(L)$ for $r,s,r_1,s_i\in \{+,-\}$, $1\le i\le 4$  and suppose further that  $e,e_1,e_4\in f$, $a_4\in yA_2z$ and $a_1\in zA_3A_1x$. A new virtual link $L_{\{e;z\}}$ called a {\it replacement} of $e$ surrounding $z$ is obtained from $L$.
Here, $L_{\{e;z\}})$
  comes from $L$ by adding vertices $x_i$ for $1\le i\le 4$ and then replacing $e$ and $e_j$ with their subdivisions for $1\le j\le 4$. If $e_{i}=e_{j}$, then $(z^{(-1)^ir_1},x_i^{-\epsilon_i})$, $(x_i^{-\epsilon_i},x_{j}^{-\epsilon_{j}})$, $(x_{j}^{-\epsilon_{j}},(-1)^{j}z^{r_1})$ is a subdivision of $e_i$ for some $1\le i\ne j\le 4$; otherwise $(z^{(-1)^ir_1},x_i^{-\epsilon_i})$, $(x_i^{-\epsilon_i},a_{i}^{s_{i}})$  is a subdivision of $e_i$.
\vskip 3mm
\noindent{\bf Lemma  $5.3.$} {\it Let $L_{\{e;z\}}$ be a {\it replacement} of $e$ surrounding $z$ as above. If either $x_i\in V(\Sigma L)$ or $z, x_i\in V(\Gamma L)$ and $\epsilon=\epsilon_i$ for $1\le i\le 4$, then $$L_{\{e;z\}}\sim L.$$}

\noindent{\bf Proof.} Because $e_i\in f$ for $0\le i\le 1$, by $\mho 0$
\begin{equation}
L_1 \sim L
\end{equation}
\noindent where $e_{x_1}=e, e_{x_3}=e_1$ and $L_1=L^{+x_1^{\epsilon_1}x_3^{\epsilon_3}}$. Here, either $x_j\in V(\Gamma L_1)$ or $x_j\in V(\Sigma L_1)$ for $j=1,3$. If $x_1,x_3\in V(\Gamma L_1)$, then $\epsilon_1=\epsilon_3.$

Because $e_{x_2}=(x_3^{\epsilon_3},y^{s})$ and $e_{x_4}=(a_4^{s_4},z^{r_1})$  are on the same face of $L_1$, by applying $\mho 0$
\begin{equation}
L_2 \sim L_1
\end{equation}
\noindent where $L_2=L_1^{+x_2^{\epsilon_2}x_4^{\epsilon_4}}$. If $x_j\in V(\Sigma L_1)$ for $j=1,3$, then $x_j\in V(\Sigma L_2)$ for $j=2,4$; otherwise $\epsilon_2=\epsilon_4$, $x_j\in V(\Gamma L_2)$ for $j=2,4$.

Because the triangle $((x_3^{\epsilon_3},x_2^{\epsilon_2}),(x_2^{-\epsilon_2},z^{r_1}),(z^{-r_1},x_3^{-\epsilon_3}))$ is an interior face of $L_2$, if either $x_i,z\in V(\Sigma L)$ for $1\le i\le 4$ or $x_i,z\in V(\Gamma L)$ and $\epsilon=\epsilon_i$ for $1\le i\le 4$,  by $\mho3$
\begin{equation}
L_2^\Delta \sim L_2
\end{equation}
where $\triangle=((x_3^{\epsilon_3},x_2^{\epsilon_2}),(x_2^{-\epsilon_2},z^{r_1}),(z^{-r_1},x_3^{-\epsilon_3}))$. If $x_i\in V(\Sigma L)$ for $1\le i\le 4$ and $z\in V(\Gamma L)$, by  $\mho4$
\begin{equation}
L_2^\Lambda \sim L_2
\end{equation}
where $\Lambda=((x_3^{\epsilon_3},x_2^{\epsilon_2}),(x_2^{-\epsilon_2},z^{r_1}),(z^{-r_1},x_3^{-\epsilon_3}))$.

 It is obvious that
\begin{equation}
L_2^\Delta =L_{\{e;z\}}\mbox{ and }L_2^\Lambda=L_{\{e;z\}}
\end{equation}
Combining (16-20), one obtains the conclusion.\hskip 10mm $\Box$
\vskip 3mm

%\vskip 2mm
Generally, let $U=\{x_i,x_i^{-}| 1\le i\le k\}$ for a nonnegative integer $k$. Given a virtual link $L$,  suppose that  $P_1=x^rx_1^{\lambda_1}x_2^{\lambda_2}\cdots x_k^{\lambda_k}y^s$ with $x\neq y \notin U$, $r,s,\lambda_i\in \{+,-\}$ for $1\le i\le k$. If either $x_i\in V(\Sigma L)$  or $x_i\in V(\Gamma L)$ and $\lambda=\lambda_i$ for  $1\le i\le k$, then $P_1$ is  called a {\it pass } of length $k$ between $x^r$ and $y^s$ in $L$.
 Add a pass $P=x^ry_1^{\delta_1}y_2^{\delta_2}\cdots y_m^{\delta_m}y^s$ on $L$ with $\delta_j\in \{+,-\}$ for $1\le j\le m$ where $y_j\notin V(L)$ for $1\le j\le m$. i.e. delete edges $e_{y_j}=(a_j^{r_j},b_j^{s_j})\in E(L)$ and then add vertices $y_j$, edges $(y_{i-1}^{\delta_{i-1}},y_{i}^{\delta_i})$, $e_{a_j}=(y_j^{-\delta_j},a_j^{r_j})$, $e_{b_j}=(y_j^{-\delta_j},b_j^{s_j})$ with $\delta_{j},r_{j},s_j\in \{+,-\}$ for $1\le j\le m$ such that $x,e_{y_1}\in f_1$, $y,e_{y_m}\in f_{m+1}$ and $e_{y_{i-1}},e_{y_{i}}\in f_i$ where $f_i$ are faces of $L$ for $1\le i\le m+1$. Here, if $e_{y_p}\ne e_{y_q}$ for any $1\le p<q\le m$, then $e_{a_j},e_{b_j}$ is the subdivision of  $e_{y_j}$ where $y_0^{\delta_{0}}=x^r$, $y_{m+1}^{\delta_{m+1}}=y^s$. Otherwise, if $e_{y_{i_1}}=e_{y_{i_2}}=\cdots=e_{y_{i_l}}$ for $l\ge 2$ and $1\le i_1<i_2<\cdots<i_l\le m$, then $(a_{i_1}^{r_{i_1}},y_{i_1}^{-\delta_{i_1}}),(y_{i_1}^{-\delta_{i_1}}, y_{i_2}^{-\delta_{i_2}}), \cdots, (y_{i_l}^{-\delta_{i_l}},b_{i_1}^{s_{i_1}})$ is the subdivision of $e_{y_{i_1}}$. Thus a planar embedding $\mu(G)$ is constructed.

\vskip 3mm

\noindent{\bf Definition $5.4.$} Let $P_1$ and $P$ be given above. If either $y_j\in V(\Gamma G)$, $\delta_{j}=\lambda_i$ for $x_i\in V(\Gamma L)$ and either there does not exist any crossing or each crossing is classical in the interior of $C$ or $y_j\in V(\Sigma G)$ for $x_i\in V(\Sigma L)$, $1\le i\le k$ and $1\le j\le m$, then a new virtual link
$L_P$
is called a {\it pass replacement } where
$C=(Px_k^{\lambda_k}x_{k-1}^{\lambda_{k-1}}\cdots x_1^{\lambda_1})$. Here,
$L_P$ forms from $\mu(G)$ by deleting vertices $x_i$ and their incident edges, and adding edges $(pre^p(x_i^{-\lambda_i}),suc^q(x_i^{-\lambda_i}))$ with $pre^p(x_i^{-\lambda_i}),suc^q(x_i^{-\lambda_i})\notin U$, $pre^l(x_i^{-\lambda_i}),suc^t(x_i^{-\lambda_i})\in U$ for each $1\le l\le p-1$, $1\le t\le q-1$ and $1\le i\le k$. Especially, if a curve consists of $x_{i_1}^{-\epsilon}, x_{i_2}^{-\epsilon},\cdots, x_{i_l}^{-\epsilon}$ for $l\ge 2$ and $1\le i_1<i_2<\cdots<i_l\le k$, then add a trivial component $O_{i_1}$ in $L_P$.
\vskip 3mm
\noindent{\bf Lemma $5.5.$} {\it For any $L\in \bar{\cal L}$, let $L_P$ be a pass replacement for any nonnegative integer $k$ and  positive integer $m$ as Definition $5.4$. Then
$$L_P\sim L.$$}
\noindent{\bf Proof.} This conclusion can be verified by using a similar argument in Lemma $3.6$ by induction on the number of crossings in the interior of $C$ and by replacing $\Omega1$ and Lemma $3.4$ with $\mho1$ and Lemma $5.3$ respectively. \hskip 10mm $\Box$
\vskip 3mm
The following conclusion is easily to deduced from $k=0$.
 \vskip 3mm

 \noindent{\bf Corollary $5.6.$} {\it Given $L\in \bar{\cal L}$, set a pass $P_1=x^{r}y_1^{\delta_1}y_2^{\delta_i}\cdots y_m^{\delta_m}y^{s}$ for $m\ge 1$ and $r,\delta_i,s\in \{+,-\}$. If $x$ and $y$ is on the same face, then
a new virtual link
$$L_P\sim L$$
where $P=x^ry^s$.}
\vskip 3mm
\noindent {\bf Lemma $5.7.$} {\it Let $f_1$ be an interior face of a connected virtual link $L$. If $f_1$ is adjacent to $f_L$, then a new virtual link $$L_1\sim L$$
 where ${\cal P}_{L_1}={\cal P}_{L}$  and $f_{L_1}=f_1.$}
 \vskip 3mm

\noindent {\bf Proof.} Using a similar argument in the proof of Lemma $3.8$ and simultaneously replacing Lemma $3.6$ with Lemma $5.5$, replacing $\Omega i$ with $\mho i$ for $0\le i\le 2$ and then  replacing $\Omega3$ with $\mho3$ or $\mho4$, one gets this case.
 \hskip 10mm $\Box$

 \vskip 3mm
Similar as  Lemma $3.9$, the following result holds.
 \vskip 3mm
\noindent{\bf Lemma $5.8.$} {\it Let $L$ be a connected virtual link. For any interior face $g$ of $L$, then there is a new virtual link $$L_g\sim L$$
such that ${\cal P}_{L_g}={\cal P}_{L}$  and $f_{L_g}=g.$}
 \vskip 3mm

This induces the following case.

\vskip 3mm
\noindent{\bf Corollay $5.9.$} {\it Cases $2$ of $\mho i$ hold for $1\le i\le 4$.}

\vskip 3mm
For a connected virtual link $L$, let $L_g$ be a virtual link in Lemma $5.8$. By Lemma $5.8$ $L_g\sim L.$ This means that $L$ and $L_g$ are virtual diagrams of the same virtual link in $S\times I$. From this viewpoint, $L_g=L$ which induces the result.

\vskip 3mm
\noindent{\bf Theorem $5.10.$} {\it For a connected virtual link $L$, $L$ is uniquely determined by ${\cal P}_{L}.$ }
\vskip 3mm

Given a disconnected virtual link $L$, similarly, $L=\bigcup\limits_{1\le i\le l}L_i$ such that $L_i$ is the shadow of $L_i$ and $f_L$ is the union of $f_{L_i}$ for $l\ge 2$ and $1\le i\le l$.  Therefore, the following result holds.

 \vskip 3mm
\noindent{\bf Theorem $5.11.$} {\it For a virtual link $L$, $L$ is uniquely determined by ${\cal P}_{L}.$ }
\vskip 3mm

 Thus, Definition $5.1$, $\mho3$ and $\mho4$, still use the same notation, are reduced as follows.
\vskip 3mm

\noindent{\bf Definition $5.1.$} {\it For two virtual links $L_i$ with $1\le i\le 2$, $L_1\cong L_2$ if and only if there exists a bijection $\phi:V(\Gamma L_1)\rightarrow V(\Gamma L_2)$ and $\phi:\{\alpha|\alpha\in V(\Sigma L_1)\}\rightarrow \{\beta^r|\beta\in V(\Sigma L_2), \mbox{ for some }r\in \{+,-\}\}$ such that $\phi(x^-)=(\phi(x))^-$ and $\phi(\sigma_{x(L_1)})=\sigma_{\phi(x)(L_2)}$ for $x\in V(L_1)$.}

\vskip 3mm

\noindent{\bf $\mho3$:} Let a triangle $\triangle=((x^{s_1},y^{s_2})$, $(y^{-s_2},z^{-s_3})$, $(z^{s_3},x^{-s_1}))$ be a face of $L$ with $s_i\in \{+,-\}$ for $1\le i\le 3$.
If either $x,y,z\in V(\Gamma L)$ for $s_1=s_2$ or $x,y,z\in V(\Sigma L)$, then
$$L^\triangle\sim L.$$
Here, ${\cal P}_{L^\triangle}={\cal P}_{L}^{\{U;W\}}$ where $U=\{x^{s_1},y^{s_2},y^{-s_2},$ $z^{-s_3},z^{s_3},$ $x^{-s_1}\}$,
$W=\{y^{s_2},x^{s_1},z^{-s_3},$ $y^{-s_2},x^{-s_1},z^{s_3}\}$.

\vskip 2mm

\noindent{\bf $\mho4$:}
If a triangle $\Lambda=((\alpha^{r},\beta^{s})$, $(\beta^{-s},a^-)$, $(a,\alpha^{-r}))$ is a face of $L$ with $r,s\in \{+,-\}$, $\alpha,\beta\in V(\Sigma L)$ and $a\in V(\Gamma L)$, then
$$L^\Lambda\sim L$$
Here, ${\cal P}_{L^\Lambda}={\cal P}_{L}^{\{U;W\}}$  where $U=\{\alpha^{r},\beta^{s},\beta^{-s},a^-,a,$ $\alpha^{-r}\}$ and
$W=\{\beta^{s},\alpha^{r},a^-,\beta^{-s},$ $\alpha^{-r},a\}$.

\vskip 3mm

Next we construct an algorithm to reduce the number of crossings in a similar way. Given a connected virtual link $L$, if either $x_i\in V(\Gamma L)$ and $r=s=-s_i\in\{+,-\}$ or  $x_i\in V(\Sigma L)$ and $x,y\in \Gamma(L)$ for $k\ge 2$ and $1\le i\le k$, then $P=x^rx_1^{s_1}\cdots x_{k}^{s_k}y^s$ is a {\it maximal pass}.
Set $G=(V,E)$ where $V(G)=V(L)\setminus \{x_i|1\le i\le k\}$ and
$$E(G)=(E(L)\cup E_1)\setminus \{(x^r,x_{1}^{-s_1}),(x_{k}^{-s_k},y^s),(x_i^{-s_i},x_{i+1}^{-s_{i+1}})|1\le i\le k-1\}.$$
Here

\hskip 10mm $E_1=\{(pre^p(x_i^{s_i}),suc^q(x_i^{s_i}))|pre^p(x_i^{s_i}),suc^q(x_i^{s_i})\notin U, pre^l(x_i^{s_i}),suc^t(x_i^{s_i})\in U $

\hskip 54mm for each $1\le l\le p-1, 1\le t\le q-1, 1\le i\le k\}$.

Let $f_x$, $f_y$ and $f_i$ denote faces obtained by removing the edges $(x^r,x_1^{-s_1})$, $(x_k^{-s_k},y^s)$ and $(x_i^{-s_i},x_{i+1}^{-s_{i+1}})$ for $1\le i\le k-1$, respectively. The dual graph of $G$ is called an {\it adjacent graph} $G_{x^ry^s}$ of $x^ry^s$.

\vskip 3mm
\noindent{\bf Reduction Crossing Algorithm I}
\vskip 2mm
Set $L$ to be  a connected virtual link with $n$ crossings, $N=0$.
\vskip 1mm
1. If there exists $e=(x,x^-)\in E(L)$ and $y^s\ne x^-$ for $s,t\in \{+,-\}$ where $(x,y^s)$ and $(x^-,z^t)$ are other adjacent edges of $x$, then let $n=n-1$, $L=L^{-x}$; otherwise go to step 2.
\vskip 1mm
2. If $(x^r,y^s),(x^{-r},y^{-s})\in E(L)$ for either $x,y\in V(\Gamma L)$ and $r=s$  or $x,y\in V(\Sigma L)$ for $r,s\in \{+,-\}$, then set $n=n-2$, $L=L^{-\{x,y\}}$ and go to step 1; otherwise go to step 3.

\vskip 1mm
3. Suppose that $P=x^rx_1^{s_1}\cdots x_{k}^{s_k}y^s$ is a maximal pass for $r,s,s_i\in \{+,-\}$ and $1\le i\le k$. Construct {\it adjacent graph} $G_{x^ry^s}$ of $x^ry^s$ and determine the shortest path $G_{x^ry^s}$ from $f_x$ to $f_y$ by Thorup's shortest algorithm. Set the length of $G_{x^ry^s}$ to be $m$. If $m\le k-1$, then let $n=n+m-k$ and $L=L_{P_{xy}}$ and go to step 1; otherwise step 5;
\vskip 1mm
4. If $m= k$ and $P_{N}\ne P_i$ for $1\le i\le N-1$, %f_xf_{x_1}\cdots f_{x_{k-1}}f_y
then let $L=L_{P_N}$, let $N=N+1$ and then go to step 1; otherwise step 5;
\vskip 1mm
5. (1) If there exists two incident edges between   $x$ and $x^-$ for each $x\in V(L)$,  then stop and $L$ is an unlink;

(2) If $L$ is an alternate link, then stop and $L$ is alternate;

(3) Otherwise, if $N\le N_0$, then let $N=N+1$ and go to step 1, otherwise stop.

(4) If $L$ is disconnected, then $L$ is a splitting link.

\vskip 5mm
 \hskip 35mm{\bf $6.$ Orientable links and virtual links }
\vskip 5mm
Given an oriented diagram  $\vec{L}$ of  an oriented link with $m$ components in $R^3$ for $m\ge 1$, label every crossing with two mutual letters $a$ and $a^-$ where $a$ and $a^-$ represent the overcrossing and undercrossing respectively.  Here, its {\it shadow}  is a marked 4-regular planar embedding  of a directed graph obtained from $\vec{L}$ by regarding $a$ and $a^-$ as the same vertex $a$(no overcrossings or undercrossings at crossings). Here, its incident arcs at $a$ are $e_i=<a^r,x_i^{r_i}>$ where there is not any other crossing between $a^r$ and $x_i^{r_i}$ along the corresponding curves of $L$ for $r,r_i\in \{+,-\}$ with $1\le i\le 4$ and $e_1e_3,e_3e_1\notin \sigma_a$.  The marked directed embedding is called {\it embedding} presentation of an oriented link. Through this paper an oriented diagram  always means an embedding presentation.

\vskip 3mm

\noindent{\bf Definition $6.1.$} For two oriented links $\vec{L}_i$ with $1\le i\le 2$, $\vec{L}_1\cong \vec{L}_2$ if and only if there exists a bijection $\phi:V(\vec{L}_1)\rightarrow V(\vec{L}_2)$ such that $\phi(a^-)=(\phi(a))^-$ and $\phi(\sigma_{a(\vec{L}_1)})=\sigma_{\phi(a)(\vec{L}_2)}$ for every $a\in V(\vec{L}_1)$.
\vskip 3mm

\vskip 3mm

%\vskip 2mm

\setlength{\unitlength}{0.97mm}
\begin{center}
\begin{picture}(100,75)

\put(28,60){\line(1,0){2}}
\put(28,58){\line(0,1){2}}
\qbezier(15,50)(-1,63)(25,61)
\qbezier(13,60)(16,42)(28,50)
\qbezier(13.5,62)(25,80)(27,58)%\qbezier(13.5,22)(25,40)(27,18)
\qbezier(17,50)(22,50)(27,58)
\qbezier(28,50)(32,56)(28,60)
\begin{footnotesize}
\put(10.5,61){{$a$}}
\put(27,61){{$b$}}
\put(18.5,48.5){{$c$}}
\put(18.5,41.5){{$\vec{L}_1$}}
\end{footnotesize}

\put(78,60){\line(1,0){2}}
\put(78,58){\line(0,1){2}}
\qbezier(67,50)(47.5,64)(77,61)
\qbezier(63,61)(66,42)(78,50)
\qbezier(63,61)(75,80)(77,58)
\qbezier(67,50)(72,50)(77,58)
\qbezier(78,50)(82,56)(77,61)
\begin{footnotesize}
\put(60.5,61){{$a$}}
\put(77,61){{$b$}}
\put(68.5,48.5){{$c$}}
\put(68.5,41.5){{$\vec{L}_1$}}
\end{footnotesize}

\put(8,20){\line(1,0){2}}
\put(10,18){\line(0,1){2}}
\qbezier(15,10)(-1,23)(25,21)
\qbezier(13,20)(16,2)(28,10)
\qbezier(13.5,22)(25,40)(27,18)
\qbezier(17,10)(22,10)(27,18)
\qbezier(28,10)(32,16)(28,20)
\begin{footnotesize}
\put(10.5,21){{$a$}}
\put(27,21){{$b$}}
\put(18.5,8.5){{$c$}}
\put(18.5,1.5){{$\vec{L}_2$}}
\end{footnotesize}

\put(58,20){\line(1,0){2}}
\put(60,18){\line(0,1){2}}
\qbezier(67,10)(47.5,24)(77,21)
\qbezier(63,21)(66,2)(78,10)
\qbezier(63,21)(75,40)(77,18)
\qbezier(67,10)(72,10)(77,18)
\qbezier(78,10)(82,16)(77,21)
\begin{footnotesize}
\put(60.5,21){{$a$}}
\put(77,21){{$b$}}
\put(68.5,8.5){{$c$}}
\put(68.5,1.5){{$\vec{L}_2$}}
\put(30,-5){{Fig.$4$: Oriented trefoils }}
\end{footnotesize}
\end{picture}
\end{center}
\vskip 8mm
\noindent{\bf Example $6.2$.} Two oriented trefoils.
\vskip 3mm
Fig.$4$ shows two oriented trefoils $\vec{L}_i$ for $1\le i\le 2$. Here,
\begin{eqnarray*}
{\cal P}_{\vec{L}_1} & = & (<b^-,a>,<b,a^->,<a,c^->,<a^-,c>)\\
              &   & (<b^-,a>,<c^-,b>,<c,b^->,<b,a^->)\\
              &  & (<c^-,b>,<a^-,c>,<a,c^->,<c,b^->),
\end{eqnarray*}
\begin{eqnarray*}
{\cal P}_{\vec{L}_2} & = & (<a^-,b>,<c^-,a>,<c,a^->,<a,b^->)\\
              &   & (<b^-,c>,<a^-,b>,<a,b^->,<b,c^->)\\
              &  & (<c^-,a>,<b^-,c>,<b,c^->,<c,a^->)
\end{eqnarray*}
%$f_{\mu(\vec{L}_1)}=(<a,c^->,<c,b^->,<b,a^->)$ and $f_{\mu(\vec{L}_2)}=(<c^-,a>,<a^-,b>,<b^-,c>)$.
\hskip 130mm $\Box$
\vskip 3mm
Let $\vec{\cal L}$ be the set of oriented links.
Now we construct an algebraic system $(\vec{\cal L},\sim)$ where an equivalence relation $"\sim"$ is defined blow.
\vskip 3mm
{\bf $\vec\Omega0$:} %$\Omega2^+$
Suppose that $e_x=<x_1^{r_1},y_1^{s_1}>$  and $e_{y}=<x_2^{r_2},y_2^{s_2}>$ are on the same face of $\vec L$ with $r_j,s_j\in \{+,-\}$ for $1\le j\le 2$. A new link  for $s\in \{+,-\}$
$$\vec L^{+x^sy^s}\sim \vec L$$
where $\vec L^{+x^sy^s}$ follows from $\vec L$ by adding vertices $x,y$ and then replacing arcs $e_x$ and $e_y$ with their subdivisions.
Here, if $e_x\ne e_y$, then $<x_1^{r_1},x^s>$, $<x^s,y^s>$, $<y^s,y_1^{s_1}>$ is the subdivision of $e_x$ and $<y_2^{s_2},x^{-s}>$, $<x^{-s},y^{-s}>$, $<y^{-s},x_2^{r_2}>$ is the subdivision of $e_y$. Otherwise,  $<x_1^{r_1},x^s>$, $<x^s,y^s>$, $<y^s,y^{-s}>$, $<y^{-s},x^{-s}>$, $<x^{-s},y_1^{s_1}>$ is the subdivision of $e_x$.
\vskip 3mm
$\vec\Omega1$: %$\Omega2^-$
Set $U=\{x,y\}$ and
set $<x^r,x^{-r}>,<y^{r},y^{-r}>\notin E(\vec L)$ for $\vec L\in \vec{\cal L}$ and $r\in \{+,-\}$.  If $<x,y>\in E(\vec L)$ and either $<x^-,y^->\in E(\vec L)$ or $<y^-,x^->\in E(\vec L)$, then
$$ \vec L^{-U}\sim \vec L.$$
Suppose that other adjacent vertices of $x^r$ and $y^r$ are $z^{s_r}$ and $v^{t_r}$ for $r\in \{+,-\}$ respectively.
$\vec L^{-U}$ is obtained from $\vec L$ by deleting vertices $x$, $y$ and their incident arcs, adding an arc $<z^{s_+},v^{t_+}>$ with $z^{s_+}\neq y$, adding an arc $<z^{s_-},v^{t_-}>$ for $z^{s_-}\neq y^-$ and $<x^-,y^->\in E(\vec L)$, adding an arc $<v^{t_-},z^{s_-}>$ for $z^{s_-}\neq y^-$ and $<y^-,x^->\in E(\vec L)$, adding two arcs $<x,x^->_i$ with $z^{s_+}= y$ and then adding two arcs $<y,y^->_i$ with $z^{s_-}= y^{-}$ for $1\le i\le 2$.

\vskip 2mm

{\bf $\vec\Omega2$:} %$\Omega1$
 Given $L\in {\cal L}$, set $<x^r,x^{-r}>\in E(L)$ for some $r\in\{+,-\}$. Suppose that other incident arcs of $x$ are $<x^r,y_r^{s_r}>$ and $<x^{-r},y_{-r}^{-r}>$. If $y_r^{s_r}\ne x^{-r}$, then
$$\vec L^{-x}\sim \vec L$$
where $\vec L^{-x}$ follows from $\vec L$ by deleting the crossing $x$ and its incident edges and then adding the edge $<y_r^{s_r},y_{-r}^{s_{-r}})$.
\vskip 2mm

{\bf $\vec\Omega3$:}
If a triangle $\triangle=(<x^{r},y^{r}>$, $<y^{-r},z^{-s}>$(or $<z^{-s},y^{-r}>)$, $<z^{s},x^{-r}>$(or $<x^{-r},z^s>$)) is a face of a link $L$ where $r,s\in \{+,-\}$, then $$L^\triangle\sim L$$
where  ${\cal P}_{L^\triangle}={\cal P}_{L}^{\{U;W\}}$ and $f_{ L^\triangle}=f_{ L}^{\{U;W\}}$. Here, $U=\{x^{r},y^{r},y^{-r},z^{-s},z^{s},x^{-r}\}$ and
$W=\{y^{r},x^{r},z^{-s},$

\noindent $y^{-r},x^{-r},z^{s}\}$.
\vskip 3mm
The following results are concluded by applying a similar method in Section $3$.
\vskip 3mm
\noindent{\bf Theorem $6.3.$} {\it For an oriented link $\vec L$, $\vec L$ is uniquely determined by ${\cal P}_{\vec L}.$ }
\vskip 3mm
\noindent{\bf Theorem $6.4.$} {\it  For any $\vec L_1,\vec L_2\in {\cal L}$,  $\vec L_1\sim \vec L_2$ if and only if $\vec L_1$ can transform into $\vec L_2$ by a sequence of $\vec\Omega i$ for $0\le i\le 3$.}
\vskip 3mm

Similarly, given an oriented virtual diagram  $\vec{L}$ of  an oriented virtual link with $m$ components in $R^3$ for $m\ge 1$, label every crossing with two mutual letters $a$ and $a^-$ where if $a\in V(\Gamma\vec{L})$, then $a$ and $a^-$ represent the overcrossing and undercrossing respectively, otherwise $a$ and $a^-$ represent two occurrences of a crossing.  Here, its {\it shadow}  is a planar embedding $\vec{L}$ of a marked $4$-regular directed graph obtained from $\vec{L}$ by regarding $a$ and $a^-$ as the same vertex $a$(no overcrossings or undercrossings at crossings) and by denoting incident arcs at $a$ with ordered pairs $<pre(a),a>,<a,suc(a)>,$ $<pre(a^-),a^->$ and $<a^-,suc(a^-)>$. Thus the {\it embedding} presentation of an oriented virtual link is obtained.

\vskip 3mm

\noindent{\bf Definition $6.5.$} For two oriented virtual links $\vec{L}_i$ for $1\le i\le 2$, $\vec{L}_1\cong \vec{L}_2$ if and only if there exists a bijection $\phi:V(\Gamma  \vec{L}_1)\rightarrow V(\Gamma \vec{L}_2)$ and $\phi:\{a|a\in V(\Sigma  \vec{L}_1)\}\rightarrow \{b^{r_b}|b\in V(\Sigma \vec{L}_2)\mbox{ and some }r_b\in \{+,-\}\}$ such that $\phi(a^-)=(\phi(a))^-$ and $\phi(\sigma_{a(\vec{L}_1)})=\sigma_{\phi(a)(\vec{L}_2)}$ for every $a\in V(\vec{L}_1)$.
\vskip 3mm
 The corresponding algebraic systems and results on oriented virtual links can directly induced  by considering the orientations on virtual links.

 \newpage
%\vskip 5mm
 \hskip 55mm{\bf $7.$ Further study }
\vskip 5mm
%\vskip 2mm
\noindent{\bf Conjecture $7.1.$ }{\it Given a $2$-connected link $L$ with $n$ crossings for $n\ge 3$, if there do not exist any short and equal pass replacement on $L$, then the crossing number of $L$ is $n$.}

\vskip 3mm
If $L$ is a alternating link with $n$ crossings,  then this conjecture holds because the Tait's Conjecture \cite{Ta98} about the crossing number of alternating link was proved independently by Kauffman \cite{Ka87}, Murasugi \cite{Mu87}, and Thistlethwaite \cite{Th87}.
Moreover, the conjecture holds for $10_{152}$   in  Rolfsen’s  tables \cite{Ro76} which is such a non-alternating knot with the least crossing number.
\vskip 3mm
\noindent{\bf Conjecture $7.2.$ }{\it For a $2$-connected link $L$ with $n$ crossings ($n\ge 3$), let $${\cal L}^*=\{L^N|L^N \mbox{ is obtained by  any }N \mbox{ equal pass replacements  based on }L, N\ge 1\}.$$  If there do not exist any $\Omega1$, $\Omega2$ and short pass replacement on $L$ for each $L^N\in {\cal L}^*$, then the crossing number of $L$ is $n$.}
\vskip 3mm

\noindent{\bf Problem $7.3.$ }{\it Let $\cal U$ be the same set of unknots as Theorem $4.1$. Does $U\in {\cal U}$ hold for any unknot $U$?}

If it is an affirmative for any unknot, then this implies that P=NP. Otherwise the problem is as follows.
\vskip 3mm
\noindent{\bf Problem $7.4.$ }{\it  Set ${\cal W}_n=\{K||V(K)|=n\}$ for $1\le n\le 2$. Let $c_1$ be a constant and  for $n\ge 3$ let
$$\begin{array}{ll}
{\cal W}_n=\{K|\mbox{ There exists }O(n^{c_1})\mbox{ path replacements such that }K\sim K_1\\
 \hskip 15mm\mbox{ and then there exists  one of }\Omega i\mbox{ and a short  path  replacement such that } \\
 \hskip 15mm K_1\sim K_2\in {\cal U}_l \mbox{ for }|\mu(K)|=|\mu(K_1)|=n, 1\le i\le 2\mbox{ and some }l<n\}.
 \end{array}
 $$
 Is there some $k\ge 1$ such that $U\in {\cal W}_k$ for any unknot $U$?}

\vskip 3mm
\noindent{\bf Problem $7.5.$} {\it Study new characteristics of a link and a virtual link by applying their embedding presentations.}
%\vskip 2mm

\vskip 3mm
\noindent{\bf Problem $7.6.$} {\it Find new presentations of a general $3$-manifold and $4$-manifold.}

If a knot $K$ is regarded as a closed non-self-intersection curve in $S^3$, then  its complement is a $3$-manifold $M_K$. Gordon and Luecke proved that a knot is determined by its complement \cite{GL89}, which means that  a knot $K$ determines the corresponding $3$-manifold $M_K$. The embedding presentation of a knot is also a presentation of $M_K$  from this point of view. How to describe a general $3$-manifold by applying the embedding presentation of a link? Furthermore, is it feasible to generalize a new presentation to a $4$-manifold?
\vskip 5mm
\noindent{\bf Acknowledgements. }
\vskip 5mm
The author is grateful to Ciprian Manolescu who provided the author an opportunity to visit Department of Mathematics at University of California, Los Angeles from Dec., 2014 to Dec., 2015. He advised the author to read some books and papers in the field of knot theory. The author would like to thank Yanpei Liu for advising the author provide an infinite family of unknots as examples. The author thank Yanxun Chang for asking the uniqueness of this representation.

This work is partially supported by NNSFC under Grant
No. 11201024 and the Foundation Grant No. 2011RC024.
\vskip 5mm
 \noindent{\bf Appendix. Unknot Haken's unknot (image courtesy of Cameron Gordon) }
\vskip 5mm
%%\begin{figure}[htp]
%\centering
%\centering{}\includegraphics[width=6in]{UnHakenU.eps}%\caption{Caption here}\label{fig1}
%\end{figure}

%\vskip 5mm

\end{document}